\documentclass[review]{elsarticle}

\usepackage[a4paper, margin=2.5cm]{geometry} 

\usepackage{amsmath,amssymb,bm}
\usepackage{amsthm}
\newtheorem{remark}{Remark}
\usepackage{subfigure}

\usepackage{xspace}
\newcommand*{\eg}{e.g.\@\xspace}

\makeatletter
\newcommand*{\etc}{%
    \@ifnextchar{.}%
        {etc}%
        {etc.\@\xspace}%
}

\makeatother

\usepackage{algorithm}
\usepackage{algpseudocode}

\usepackage[normalem]{ulem}
\usepackage{soul}

\usepackage[dvipsnames]{xcolor}
\usepackage[colorinlistoftodos, obeyFinal]{todonotes}

\newcommand{\abs}[1]{\lvert#1\rvert}

\journal{Railway Engineering Science}

\usepackage{setspace}
\begin{document}

\setcounter{tocdepth}{1}
%

\begin{frontmatter}

	\title{A General and Robust 3D Finite Element Dynamics Framework for Railway Vehicle–Bridge Interaction with Nonlinear Wheel–Rail Contact Modeling}

\author[epfl]{Pablo Antolin\corref{cor1}}
\ead{pablo.antolin@epfl.ch}

\author[etsiae]{Khanh Nguyen}
\ead{khanhnguyen.gia@upm.es}

\author[caminos]{Jos\'e M. Goicolea}
\ead{jose.goicolea@upm.es}

\cortext[cor1]{Corresponding author.}

\address[epfl]{Institute of Mathematics,\\ \'Ecole Polytechnique F\'ed\'erale de Lausanne, CH-1015 Lausanne, Switzerland}
\address[etsiae]{ETS de Ingenier\'{i}a Aeron\'{a}utica y del Espacio, Universidad
Polit\'{e}cnica de Madrid, Ciudad Universitaria, 28040 Madrid, Spain}
\address[caminos]{ETS de Caminos, Canales y Puertos, Universidad
Polit\'{e}cnica de Madrid, Ciudad Universitaria, 28040 Madrid, Spain}

\begin{abstract}
A key challenge in 3D finite element models of coupled railway vehicle-bridge dynamics is the rigorous definition of kinematic constraints and the development of an efficient, robust solution. This paper presents a novel approach that can be implemented in general finite element software using constraint equations tailored to wheel-rail contact behavior, essential for analyzing lateral vehicle-bridge interactions. The method employs absolute coordinates to describe the motion of nodes defining the track position and orientation for each wheelset, without assuming infinitesimal displacements or rotations. This general formulation enables realistic simulations of extreme scenarios involving large lateral movements caused by strong winds or earthquakes. The proposed wheel-rail contact model is first validated against published results, and a 3D numerical example demonstrates the method's performance and capabilities.
\end{abstract}


\end{frontmatter}


\section{Introduction} \label{sec:introduction}
The dynamic interaction between moving trains and railway bridges (Vehicle-Bridge Interaction, VBI) is a critical aspect of modern railway engineering, especially in terms of structural safety and integrity, running safety, and passenger comfort under a wide range of operating and environmental conditions. This field has garnered extensive attention over recent decades, with numerous analytical, numerical, and experimental investigations aimed at unraveling the complex dynamic effects of trains moving on bridges. The continuous push towards higher passenger train speeds and longer bridge spans has further intensified the need for accurate VBI analysis, not only for innovative railway bridge design, but also for safeguarding train operations and enhancing passenger ride quality. Consequently, the development of precise and computationally robust VBI models remains a paramount challenge for structural and mechanical engineers. Recent reviews, such as those of \cite{Zhai2019} and \cite{Montenegro2021}, continue to highlight these evolving challenges and advancements in the field.  


VBI models vary in complexity depending on the study's objectives, ranging from focusing solely on vertical dynamics to encompassing coupled lateral-torsional responses. Vehicle representations can span from simplified models, such as moving loads, moving-mass models, or sprung-mass models \cite{Fryba1996, Yang2004}, to sophisticated multibody system (MBS) models \cite{Yang2004, Tanabe2003, Xia2003, Xu2004, Shabana2008, Nguyen2009, Antolin2012, Antolin2013, Zhai2013, Qiao2018, Xiao2019, Liu2020, Montenegro2022a, Pan2024} that offer a more realistic depiction of train components such as car bodies, bogies, and wheelsets as rigid bodies with detailed inertial properties. While many traditional MBS vehicle models assume linear behavior with small displacements and rotations, the demand for simulating extreme events and complex nonlinear phenomena necessitates formulations that can accommodate large motions and geometric nonlinearities.


Bridge structures are typically discretized using the Finite Element Method (FEM), which offers detailed and flexible modeling of their dynamic behavior. Various FEM model types are described in the literature, each tailored to different applications and levels of complexity. Simplified models, such as beams and trusses, are commonly used due to their simplicity and reliability for preliminary analyzes, as demonstrated in studies like \cite{Yang2004, Tanabe2003, Zhang2008, Nguyen2009, Antolin2012, Zhai2013, Majka2008, Kwark2004}. For more detailed structural analyzes, shell elements \cite{Song2003, Oliva2013} are used to capture surface-level behaviors, while solid elements \cite{Kwasniewski2006} provide a comprehensive representation of three-dimensional stress states within the structure. In addition, hybrid models, such as beam-shell or beam-solid configurations \cite{Henchi1998}, are sometimes adopted to balance computational efficiency and modeling accuracy. These approaches offer flexibility, allowing engineers to select the appropriate level of detail based on the specific needs of their analysis. However, detailed models are costly to develop and more prone to errors, making the careful selection of a suitable model a critical consideration. 

\begin{figure}[!ht]
 \centering
 \includegraphics[width=\textwidth]{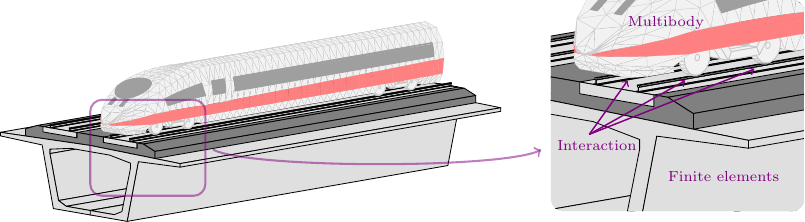}
 \caption{Sketch of the different components of a vehicle-structure model}
 \label{fig:interaction}
\end{figure}

The core of the VBI analysis lies in modeling the interaction, primarily at the wheel-rail interface and the track-bridge connection (see Fig.~\ref{fig:interaction}). The Wheel-rail interaction is especially critical for studies that involve lateral dynamics, such as derailment risks, hunting motion, and other dynamic phenomena. The existing literature categorizes wheel-rail interaction models as follows. 
\begin{itemize}
    \item Direct and simple models assume a perfectly guided path for the vehicle's wheelsets, meaning that the contact points between the wheels and rails are fixed in position and share the same velocity. Such models, as described in \cite{Xia2003,Zhang2008,Song2003,Neves2014}, are relatively straightforward. Furthermore, hunting motion, a key lateral dynamic behavior, can be introduced by prescribing the wheelset's movement \cite{Xia2003,Xu2004,Zhang2008}.
    \item Linear models incorporating relative motion between the wheelsets and the track, assuming constant conicity in the profiles of the wheels and rails. Contact mechanics are simplified, making these models computationally efficient for many applications \cite{Yang2004,Xu2004,Zhang2010}, however, this approach is restricted to analysis under normal operating conditions.
    \item Nonlinear models: These are more advanced and provide a detailed representation of the wheel, rail, and their interaction. Sophisticated techniques, such as those proposed in \cite{Tanabe2003,Nguyen2009,Antolin2012,Antolin2013,Olmos2013,Montenegro2015,Dimitrakopoulos2015,Montenegro2023}, account for the actual geometries of wheel and rail surfaces. These models enable precise localization of contact points, whether online or offline, improving the accuracy of dynamic analyzes.  
\end{itemize}

In the third approach, the accurate prediction of the location of the contact points is the key aspect to achieve any formulation for the wheel-rail contact problem. Chen and Zhai \cite{Chen2004} proposed a so-called wheel/rail spatially dynamic coupling model, in which only a contact point at each wheel/rail pair is determined by computing the minimum vertical distance between the left and right wheel/rail in real time. Pombo et al.\ \cite{Pombo2007} considered the wheel and rail surfaces as sweep surfaces obtained by dragging plane curves along spatial curves; as a result, the locations of the contact points were determined using the two-dimensional profiles of the wheel and rail. By using an optimized search technique, they could localize two possible contact points: one on the wheel tread and the other on the wheel flange. However, this methodology \cite{Pombo2007} requires that the wheel and rail surfaces be convex, and that is not always true of their actual profiles. A similar technique is also used by \cite{Montenegro2015}, but only one contact point is considered for each wheel-rail pair. Montenegro et al.\ \cite{Montenegro2023} developed an enhancement of the previous technique that simplifies the detection of contact points, particularly in the concave regions of the thread-flange transition. This enhancement enables more precise simulations of wheel-rail interactions across a range of operational conditions. However, the wheel-rail contact model remains confined to a 2D plane, limiting its ability to identify contact points at varying longitudinal positions.

The track-bridge interaction also plays a critical role in overall system dynamics. The vibration and deformation of the bridge affect the track structure through the track-bridge interface, altering the wheel-rail contact geometry and the corresponding contact forces. The track structure acts as a low-pass filter, i.e., absorbs high-frequency energy generated by wheel-rail contact. Consequently, the models that include the track structure typically predict lower acceleration amplitudes in the high-frequency range ($>$20 Hz) for the bridge deck compared to models where wheel-rail contact loads directly to the bridge \cite{Zhai2019}. The track-bridge interaction is typically modeled using finite element methods (FEM) to represent all track and bridge components, with kinematic and dynamic relationships established between them for comprehensive analysis \cite{Wu2001,Zhai2013,Olmos2013}. Such detailed modeling approaches provide deeper insight into complex train-bridge dynamics and are essential for ensuring safety and optimal performance.

Two primary strategies exist for solving the transient VBI problem. The first one, also known as the uncoupled method, treats the vehicle-bridge system as two separate subsystems. Each subsystem is governed by its uncoupled equations of motion, which are solved independently. Kinematic and dynamic compatibility between the two subsystems is enforced iteratively at each time step \cite{Xu2004,Zhang2008,Nguyen2009,Zhang2010}. The second strategy considers the vehicle-bridge system as a monolithic coupled entity, represented by a unified set of equations of motion. These equations are solved simultaneously as an integrated system \cite{Wu2001,Xia2003,Yang2004,Yau2004,Antolin2012,Antolin2013,Zhai2013,Montenegro2015,Dimitrakopoulos2015,Montenegro2023}. A detailed analysis of the advantages and limitations of these approaches is provided in \cite{Yang1996, Zhai2019}. While significant progress has been made, a persistent challenge lies in developing a robust and general kinematic constraint formulation within a monolithic, non-linear FEM framework capable of handling large 3D displacements and rotations inherent in extreme operational scenarios (e.g., strong winds, earthquakes) or complex track geometries (e.g., curved/banked bridges). Many existing approaches simplify these aspects or are implemented in specialized rather than general-purpose FE software.

This work aims to directly address this gap by proposing a novel kinematic constraint formulation for 3D vehicle-bridge interaction, which includes an accurate nonlinear wheel-rail contact model.  Our primary contribution is the development of a rigorous method using absolute coordinates to describe motion and interaction, which is implemented within general-purpose finite element software. Key features of this formulation include the following:
\begin{itemize}
    \item The introduction of a set of virtual nodes that precisely characterize the track position and orientation relative to each wheelset of the vehicle.
    \item A formulation that inherently accommodates large displacements and rotations, avoiding simplifications tied to infinitesimal motion assumptions.
    \item A 3D wheel-rail contact model is capable of accurately representing the real geometry of wheels and rails, handling multi-contact scenarios at different positions (such as the rail head, flange, or both), and accounting for variations in surface curvature within the contact zones.
    \item The ability to capture detailed track geometry, including irregularities and cant deficiency, and its deformation under load
    \item A modular design facilitates its integration into existing FE software with multibody capabilities
\end{itemize}

This approach provides a robust framework for developing VBI models capable of accurately capturing key phenomena governing lateral vehicle–structure interaction in extreme scenarios, particularly the transition between contact and separation (wheel lift) under combined lateral and vertical loading. Such capability is essential for the analysis of lateral vehicle–bridge interaction, where loss of contact and re-contact may occur as a result of structural deformation, track alignment irregularities, and dynamic amplification effects. Importantly, the proposed framework can be easily extended to cases where the track structure is explicitly included in the FEM, allowing detailed investigation of high-frequency track responses without modification of the methodology. The vehicles are modeled as three-dimensional multibody systems, and the bridge structure is discretized using finite elements such as beam, shell, or solid elements. The formulation is designed to fully consider geometrical and material nonlinearities in both subsystems. The paper is structured as follows. Section \ref{sec:interaction} will outline the general equations of motion for the coupled system. Section \ref{sec:kinematics} provides a detailed description of the proposed kinematic interaction. Section \ref{sec:contact} describes the proposed wheel-rail contact model. The proposed methodology is then applied to two different test cases in Section \ref{sec:examples}. Finally, the main conclusions are summarized in Section \ref{sec:conclusions}. Appendices detail the linearization of constraints and implementation aspects.

\section{Framework for Vehicle-Bridge Interaction} \label{sec:interaction}
For analyzing the coupled dynamic interaction between a deformable structure
and a train vehicle (see Fig.~\ref{fig:interaction}), it is necessary to formulate the equations that govern the response of both subsystems as well as their interaction.
\par 
In general terms, the equation that governs the dynamic response of a
structure modeled by means of finite elements can be expressed as:
\begin{subequations} \label{eq:structure}
\begin{align}
  \mbox{\boldmath$\sf M$}^S\ddot{\mbox{\boldmath$\sf q$}}_S
+ \mbox{\boldmath$\sf F$}^S_{\text{damp}}(\mbox{\boldmath$\sf q$}_S,\dot{\mbox{\boldmath$\sf q$}}_S)
+ \mbox{\boldmath$\sf F$}^S_{\text{def}}(\mbox{\boldmath$\sf q$}_S,\dot{\mbox{\boldmath$\sf q$}}_S)
+ \mbox{\boldmath$\sf F$}^S_{\text{SV}}(
\mbox{\boldmath$\sf q$}_S,\dot{\mbox{\boldmath$\sf q$}}_S,
\mbox{\boldmath$\sf q$}_V,\dot{\mbox{\boldmath$\sf q$}}_V)
 &= \mbox{\boldmath$\sf F$}^S_{\text{ext}}(t)\,,\\
\label{eq:structure:restriction}
\mbox{\boldmath$\sf \Phi$}^S(\mbox{\boldmath$\sf q$}_S) &= \mbox{\boldmath$\sf 0$}\,,
\end{align}
\end{subequations}
being $\mbox{\boldmath$\sf q$}_S$ the generalized coordinates vector of the structure system and $\dot{\mbox{\boldmath$\sf q$}}_S$ and $\ddot{\mbox{\boldmath$\sf q$}}_S$ its first and second time derivatives, respectively; $\mbox{\boldmath$\sf q$}_V$, $\dot{\mbox{\boldmath$\sf q$}}_V$ and $\ddot{\mbox{\boldmath$\sf q$}}_V$ are their counterparts for the
vehicle; $\mbox{\boldmath$\sf M$}^S$ is the mass matrix of the structure; 
and $\mbox{\boldmath$\sf \Phi$}^S = \mbox{\boldmath$\sf 0$}$ the set of constraint equations of the structure model.
The force vectors appearing in this equation are:
\begin{description}
\item[$\mbox{\boldmath$\sf F$}^S_{\textnormal{damp}}$:] structural
    damping force vector.
\item[$\mbox{\boldmath$\sf F$}^S_{\textnormal{def}}$:] internal forces that
    correspond to the elastic (or inelastic) deformation of the structure.
\item[$\mbox{\boldmath$\sf F$}^S_{\textnormal{SV}}$:] wheel-rail
    contact forces representing 
    the dynamic interaction between the structure and the vehicle running
    across it.
\item[$\mbox{\boldmath$\sf F$}^S_{\textnormal{ext}}$:] external forces,
    \eg, structure self-weight, aerodynamic loads, seismic actions, \etc.
\end{description}
The term \emph{force vector} is used ambiguously to refer to the vectors $\mbox{\boldmath$\sf F$}$, which include both forces and moments, as the generalized coordinates vector encompasses both displacement and rotational degrees of freedom.

\begin{remark}
It is worth noting that Eq.~\eqref{eq:structure} remains valid when the track structure is explicitly incorporated into the finite element model together with the bridge structure.
The extra degrees-of-freedom associated to the track subsystem would be part of the vector  $\mbox{\boldmath$\sf q$}_S$.
In such a case, the interaction between the vehicle and the infrastructure would be established between the track system and the vehicle subsystem, and its contribution would be incorporated in the force vector $\mbox{\boldmath$\sf F$}^S_{\textnormal{SV}}$.
\end{remark}

Similarly to the structure, the dynamic equation for a generic multibody model of a train crossing a bridge can be expressed as:
\begin{subequations}\label{eq:vehicle}
\begin{align}
  \mbox{\boldmath$\sf M$}^V\ddot{\mbox{\boldmath$\sf q$}}_V
+ \mbox{\boldmath$\sf F$}^V_{\text{susp}}(\mbox{\boldmath$\sf
q$}_V,\dot{\mbox{\boldmath$\sf q$}}_V)
+ \mbox{\boldmath$\sf F$}^V_{\text{gir}}(\mbox{\boldmath$\sf
q$}_V,\dot{\mbox{\boldmath$\sf q$}}_V)
+ \mbox{\boldmath$\sf F$}^V_{\text{SV}}(
\mbox{\boldmath$\sf q$}_V,\dot{\mbox{\boldmath$\sf q$}}_V, 
\mbox{\boldmath$\sf q$}_S,\dot{\mbox{\boldmath$\sf q$}}_S)
 &= \mbox{\boldmath$\sf F$}^V_{\text{ext}}(t),\\
\label{eq:vehicle:restriction}
\mbox{\boldmath$\sf \Phi$}^V(\mbox{\boldmath$\sf q$}_V,\dot{\mbox{\boldmath$\sf
q$}}_V,t) &= \mbox{\boldmath$\sf 0$}\,,
\end{align}
\end{subequations}
being $\mbox{\boldmath$\sf M$}^V$ the mass matrix of the vehicle model; 
and $\mbox{\boldmath$\sf \Phi$}^V = \mbox{\boldmath$\sf 0$}$ the set of constraint equations of the multibody model of the vehicle. Unlike the structural restriction Eq.~\eqref{eq:structure:restriction}, the constraint equations for the train's multibody model are formulated more generally, also depending on velocity terms and time. The force vectors in Eq.~\eqref{eq:vehicle} are:
\begin{description}
\item[$\mbox{\boldmath$\sf F$}^V_{\textnormal{susp}}$:] forces due to the
    suspension systems of the vehicle.
\item[$\mbox{\boldmath$\sf F$}^V_{\textnormal{gir}}$:] gyroscopic terms consequence of
  the rigid body rotations.
\item[$\mbox{\boldmath$\sf F$}^V_{\textnormal{SV}}$:] wheel-rail
    contact forces representing
    the dynamic interaction between the structure and the vehicle. It is
    the counterpart of $\mbox{\boldmath$\sf F$}^S_{\textnormal{SV}}$.
\item[$\mbox{\boldmath$\sf F$}^V_{\textnormal{ext}}$:] external forces,
    \eg, self-weight, aerodynamic loads, \etc.
\end{description}
\par 
Both subsystems (Eqs.~\eqref{eq:structure} and \eqref{eq:vehicle}) can be combined into a single set of equations expressed as:
\begin{subequations}\label{eq:full_system}
\begin{align}
  \mbox{\boldmath$\sf M$}\,\ddot{\mbox{\boldmath$\sf q$}}
+ \mbox{\boldmath$\sf F$}_{\text{int}}(\mbox{\boldmath$\sf q$},\,\dot{\mbox{\boldmath$\sf q$}})
+ \mbox{\boldmath$\sf F$}_{\text{SV}}(\mbox{\boldmath$\sf q$},\,\dot{\mbox{\boldmath$\sf q$}})
 &= \mbox{\boldmath$\sf F$}_{\text{ext}}(t)\,,\\
\mbox{\boldmath$\sf \Phi$}(\mbox{\boldmath$\sf q$},\,\dot{\mbox{\boldmath$\sf
q$}},\,t) &= \mbox{\boldmath$\sf 0$}\,,\\
\label{eq:full_system:restriction2}
 \mbox{\boldmath$\sf \Phi$}_{SV}(\mbox{\boldmath$\sf q$}) &= \mbox{\boldmath$\sf 0$}\,,
\end{align}
\end{subequations}
where $\mbox{\boldmath$\sf q$}$ is the generalized coordinates vector of the
full system, that can be expressed as $\mbox{\boldmath$\sf q$}=
\begin{bmatrix}\mbox{\boldmath$\sf q$}^{\text{T}}_{S},\,\mbox{\boldmath$\sf
q$}^{\text{T}}_{V}\end{bmatrix}^{\text{T}}$;
$\mbox{\boldmath$\sf M$}$ is the mass matrix;
$\mbox{\boldmath$\sf F$}_{\text{int}}$ arranges the previously defined vectors
$\mbox{\boldmath$\sf F$}^S_{\textnormal{damp}}$,
$\mbox{\boldmath$\sf F$}^S_{\textnormal{def}}$,
$\mbox{\boldmath$\sf F$}^V_{\textnormal{susp}}$ and
$\mbox{\boldmath$\sf F$}^V_{\textnormal{gir}}$;
and, in the same way,
$\mbox{\boldmath$\sf F$}_{\text{SV}}$ and
$\mbox{\boldmath$\sf F$}_{\text{ext}}$
are composed of 
$\mbox{\boldmath$\sf F$}^S_{\textnormal{SV}}$ and
$\mbox{\boldmath$\sf F$}^V_{\textnormal{SV}}$, and
$\mbox{\boldmath$\sf F$}^S_{\textnormal{ext}}$ and
$\mbox{\boldmath$\sf F$}^V_{\textnormal{ext}}$, respectively.
The restriction equation
$\mbox{\boldmath$\sf \Phi$}(\mbox{\boldmath$\sf q$},\,\dot{\mbox{\boldmath$\sf q$}},\,t) = \mbox{\boldmath$\sf 0$}$
gathers Eqs.~\eqref{eq:structure:restriction} and \eqref{eq:vehicle:restriction}; 
and 
$\mbox{\boldmath$\sf \Phi$}_{SV}(\mbox{\boldmath$\sf q$}) = \mbox{\boldmath$\sf 0$}$
correspond to the restriction equations needed for establishing the interaction between the
structure and the vehicle. The main purpose of this work is to describe the kinematic relationships that are needed to establish the set of constraint Eqs.~\eqref{eq:full_system:restriction2}, as detailed in Section \ref{sec:kinematics}.


\section{Kinematic Constraints Formulation} \label{sec:kinematics}

\begin{figure}[ht]
 \centering
\subfigure[The virtual nodes $r_1$, $r_2$ and $m$ allow to establish the kinematic coupling between 
 finite elements of the structure model and the vehicle wheelset]{\label{fig:virtual_nodes}\includegraphics[width=0.45\textwidth]{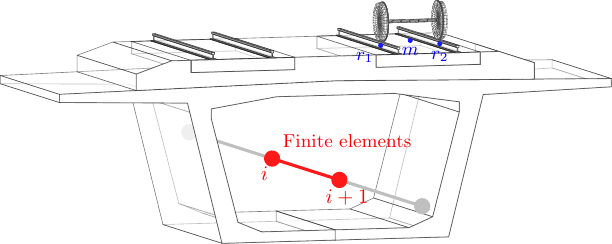}}
\subfigure[Sketch of the longitudinal positions of the nodes of the finite elements, the
 moving node $m$, and orientation vectors of the beam section]{\label{fig:track}\includegraphics[width=0.5\textwidth]{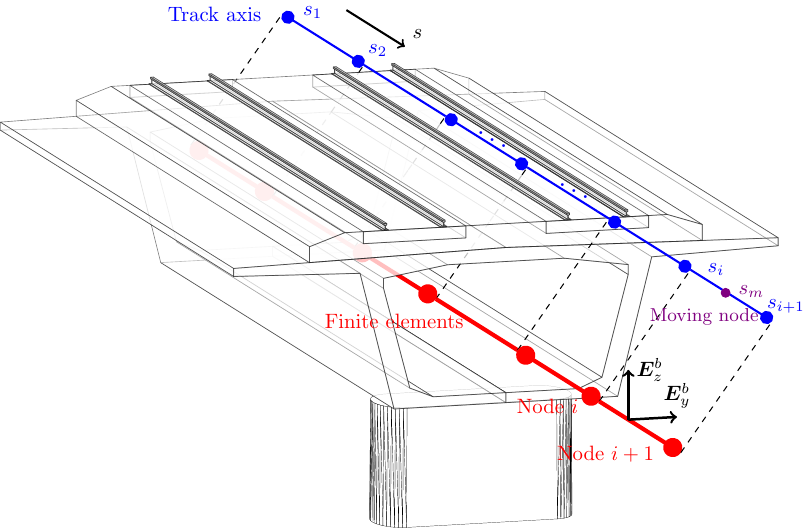}} 

 \caption{The virtual nodes $r_1$, $r_2$ and $m$ allow to establish the kinematic coupling between 
 finite elements of the structure model and the vehicle wheelset}
 
\end{figure}

Kinematics is crucial for defining the interaction between a train and a deformable structure. In this work, the bridge structure is assumed to be modeled as a 3D beam model. However, the method proposed in this work can be extended straightforwardly to models of shells, solid elements, or combinations of them. Three virtual nodes are introduced for each train wheelset to establish the kinematic relationships between the subsystems (see nodes $r_1$, $r_2$, and $m$ in Fig.~\ref{fig:virtual_nodes}). These massless nodes, existing solely for kinematic description, are positioned within the same cross-section of the wheelset and follow the trajectory defined by the track.
\par
The first virtual node, $m$, is positioned between the two rails and represents the midpoint of the track relative to the finite elements of the deformable structure. The other two virtual nodes, $r_1$ and $r_2$, are located at the centroids of each rail in the same cross-section as node $m$, enabling an accurate description of the kinematics of each rail independently, while accounting for possible track irregularities.

In the following Section \ref{sec:constraint1}, the kinematic relationship of the virtual node $m$ with respect to the finite elements of the structure is described. In the same way, in Section \ref{sec:constraint2}, the link between nodes $r_j$ and $m$ is developed. Here, the subscript $j$ refers interchangeably to $j=1$ (right rail) and $j=2$ (left rail).

\subsection{Constraint between the moving node and the structure} \label{sec:constraint1}

This section defines the kinematic relationships between node $m$ and the finite elements of the structure. The virtual node $m$, moving at the same velocity $v$ as the train, is located in the same cross-section of the structure as the wheelset, as illustrated in Fig.~\ref{fig:track}.

To achieve this, the track's longitudinal direction is parameterized using the variable $s$. Assuming a wheelset is initially at the track section defined by $s_0$, the position of node $m$ at time $t$, denoted as $s_m$, is calculated as:
\begin{align}\label{eq:sm}
    s_m = s_0+\int^t_0v\,\text{d}t\,.
\end{align}
If $v$ is constant, Eq.~\eqref{eq:sm} simplifies to $s_m = s_0 + vt$. As shown in Fig.~\ref{fig:track}, node $m$ lies on track segment $i$ (between beam nodes $i$ and $i+1$) if $s_i \leq s_m < s_{i+1}$, where $s_i$ and $s_{i+1}$ are the longitudinal positions of two consecutive structural nodes. To identify the finite element in which node $m$ is located, the longitudinal deformation of the structure is neglected, as it is insignificant relative to $s$.
\begin{figure}[ht]
 \centering
 \includegraphics[width=0.7\textwidth]{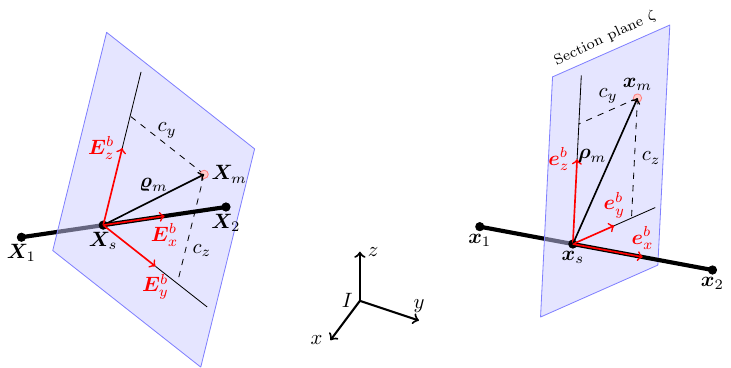}
 \caption{Orientation vectors of the beam elements and offsets of the node $m$
 for the initial (left) and current (right) configurations}
 \label{fig:beams}
\end{figure}
\par 
In the reference (or undeformed) configuration, the unit vectors $\bm{E}^b_y$ and $\bm{E}^b_z$ orient the beam section (see Figs.~\ref{fig:track} and \ref{fig:beams}). They are two unit vectors defined such that $\bm{E}^b_z$ is normal to the track plane in every section\footnote{In the case of curved and/or banked bridges, $\bm{E}^b_z$ has the same direction as the bi-normal vector of the track axis curve.} (pointing upwards); and $\bm{E}^b_y$ is defined such that the beam tangent vector $\bm{E}^b_x$ (defined in the direction of $\bm{X}_2-\bm{X}_1$) along with $\bm{E}^b_y$ and $\bm{E}^b_z$ forms a right-handed orthonormal basis. Hereinafter, all the vectors are referred to a global inertial reference system (denoted with $Ixyz$ in Fig.~\ref{fig:beams}).
\par 
The initial offsets $c^m_y$ and $c^m_z$ of the node $m$, which is initially located at a given section plane $\zeta$ respect to the beam axis, can be expressed as 
\begin{subequations} \label{ch3:eq:excen}
\begin{align}
 c^m_y & = \left(\bm{X}_m-\bm{X}_1\right)\cdot\bm{E}^{b,0}_y\,,\\
 c^m_z & = \left(\bm{X}_m-\bm{X}_1\right)\cdot\bm{E}^{b,0}_z\,,
\end{align}
\end{subequations}
where $\bm{X}_m$ is the initial position vector of node $m$, $\bm{X}_1$, which corresponds to the first node of the beam element, where the wheelset is initially located, and $\bm{E}^{b,0}_y$ and $\bm{E}^{b,0}_z$ are the orientation vectors defined above but referred to the initial position of the wheelset.
\par 
The position $\bm{X}_m$ can be expressed as (see Fig.~\ref{fig:beams})
\begin{align} \label{ch3:eq:Xm}
\bm{X}_m = \bm{X}_s + \bm{\varrho}_m\,,
\end{align}
where $\bm{X}_m$ is the point of the beam element in the cross section
$s_m$ (see Figs.~ \ref{fig:track} and \ref{fig:beams}), and $\bm{\varrho}_m$ is:
\begin{align} \label{ch3:eq:Rhom}
  \bm{\varrho}_m = c^m_y\,\bm{E}^b_y + c^m_z\,\bm{E}^b_z\,.
\end{align}
Similarly, the position of $m$ in the current configuration can be expressed as:
\begin{align}\label{eq:xm}
 \bm{x}_m = \bm{x}_s+\bm{\rho}_m\,,
\end{align}
where $\bm{x}_s$ is position of the beam element at section $s_m$ and
$\bm{\rho}_m$ is the counterpart of $\bm{\varrho}_m$, in the 
deformed configuration.
\par 
Assuming that the structural section remains undeformed in the transverse directions $\bm{E}^b_y$ and $\bm{E}^b_z$, the initial offsets $c^m_y$ and $c^m_z$ will be constant. Consequently, $\bm{\rho}_m$ can be expressed as:
\begin{align}\label{ch3:eq:rest_disp1}
 \bm{\rho}_m = c^m_y\,\bm{e}^b_y +c^m_z\,\bm{e}^b_z\,,
\end{align}
where $\bm{e}^b_y$ and $\bm{e}^b_z$ are the corresponding vectors $\bm{E}^b_y$ and $\bm{E}^b_z$, referred to the current configuration. They can be formulated as a rotation of $\bm{E}^b_y$ and $\bm{E}^b_z$
\begin{align}\label{ch3:eq:rest_rot}
  \bm{e}^b_{\alpha} = \bm{\Lambda}(\bm{\phi}_m)\bm{E}^b_{\alpha}\,,\quad \alpha=\lbrace
    y,\,z\rbrace\,,
\end{align}
being $\bm{\phi}_m$ the rotation vector of the beam section $s_m$
(introduced below) and $\bm{\Lambda}_m =\bm{\Lambda}(\bm{\phi}_m)$ its associated
rotation tensor (detailed in \ref{sec:kinematic_linearization}). In the same way, $\bm{\rho}_m$ can be expressed simply as:
\begin{align}\label{ch3:eq:rest_disp1_bis}
  \bm{\rho}_m = \bm{\Lambda}_m\bm{\varrho}_m\,.
\end{align}
\par 
The position of a point in the deformed configuration can be written as
$\bm{x} = \bm{X} + \bm{u}$, where $\bm{X}$ and $\bm{u}$ are its initial position and displacement vectors.
Thus, Eq.~\eqref{eq:xm} can be expressed as
\begin{align}\label{eq:rest_disp2_}
  \bm{x}_m = \bm{X}_m + \bm{u}_m = \bm{X}_s + \bm{u}_s + \bm{\rho}_m\,,
\end{align}
and therefore
\begin{align}\label{eq:rest_disp2_bis}
  \bm{u}_m = \bm{X}_s + \bm{u}_s + \bm{\rho}_m - \bm{X}_m\,.
\end{align}
\par 
Being $s_1$ and $s_2$ the longitudinal parameters of two nodes of a beam element on which the virtual node $m$ is located, its displacement vector can be expressed, as a function of its natural coordinate $\zeta = \left(2\,s_m-s_1-s_2\right)/\left(s_2-s_1\right)$, as\footnote{Other possible element formulations, with a different number of nodes, such as quadratic beam
elements or shells, could also be used.}
\begin{align}\label{eq:rest_disp3}
  \bm{u}_m &= N_1(\zeta) \left(\bm{X}_1 + \bm{u}_1 \right)+ N_2(\zeta) \left(\bm{X}_2 + \bm{u}_2 \right) + \bm{\rho}_m - \bm{X}_m\,,
\end{align}
where $\bm{u}_1$ and $\bm{u}_2$ are the displacement vectors of both nodes of the beam element. $N_1$ and $N_2$ are the basis functions of both nodes of the beam element evaluated at $\zeta$, and $\bm{X}_1$, $\bm{X}_2$ are their nodal positions in the initial configuration.
\par 
On the other hand, the orientation of both nodes in the deformed configuration
is expressed through their rotation vectors $\bm{\phi}_1$ and $\bm{\phi}_2$. The magnitude of a rotation vector represents the angle of rotation around its own direction, known as the Euler axis (see \cite{Shabana2001} for details). At section $s_m$, defined by the natural coordinate $\zeta$, the rotation vector $\bm{\phi}_m$ is interpolated as:
\begin{align}\label{eq:interp_rotation}
\bm{\phi}_m = N_1(\zeta)\,\bm{\phi}_1+N_2(\zeta)\,\bm{\phi}_2\,.
\end{align}
The interpolation formula \eqref{eq:interp_rotation} is neither exact nor the only possible approach. A detailed discussion of rotation interpolation methods and their implications can be found in \cite{Romero2004}. As noted in \cite{Crisfield1998}, the interpolation in \eqref{eq:interp_rotation} has two main limitations: the resulting elastic operators lack material objectivity and depend on the sequence of rotations. However, these drawbacks are not significant in this work, as the differences between $\bm{\phi}_1$ and $\bm{\phi}_2$ are small, and no elastic operators are obtained using $\bm{\phi}_m$.
\par 
Consequently, the kinematics of the node $m$ can be expressed as a function of the generalized coordinates of the beam nodes through Eqs.
\eqref{eq:rest_disp3} and \eqref{eq:interp_rotation}.
\par 
In a broader sense, the generalized coordinates of all the nodes $m$ (one
for each wheelset of the train), denoted as $\mbox{\boldmath$\sf q$}_m$, which are dependent, can be expressed as a function of the generalized coordinates of the structure nodes $\mbox{\boldmath$\sf q$}_S$, the independent coordinates, as follows:
\begin{align}\label{eq:restricciongeneral}
\mbox{\boldmath$\sf q$}_m &= \mbox{\boldmath$\sf\Xi$}^m(\mbox{\boldmath$\sf q$}_S)\,.
\end{align}
This can be rewritten as a restriction equation
$\mbox{\boldmath$\sf \Phi$}^m_{SV}= \mbox{\boldmath$\sf q$}_m
- \mbox{\boldmath$\sf\Xi$}^m(\mbox{\boldmath$\sf q$}_S)
= \mbox{\boldmath$\sf0$}$, that is included in the Equation
\eqref{eq:full_system:restriction2}.

Most of the available methods for imposing constraints require a linearization
of the restriction equations. Thus, $\mbox{\boldmath$\sf q$}_m$ is linearized
as:
\begin{align}\label{eq:restricciongeneral3}
\delta\mbox{\boldmath$\sf q$}_m &=
\frac{\partial\mbox{\boldmath$\sf\Xi$}^m(\mbox{\boldmath$\sf
q$}_S)}{\partial\mbox{\boldmath$\sf q$}_S}\,\delta\mbox{\boldmath$\sf q$}_S\,,
\end{align}
where $\delta\mbox{\boldmath$\sf q$}_m$ and $\delta\mbox{\boldmath$\sf q$}_S$ are
the variations of $\mbox{\boldmath$\sf q$}_m$ and $\mbox{\boldmath$\sf q$}_S$, and
the operator $\partial\mbox{\boldmath$\sf\Xi$}^m(\mbox{\boldmath$\sf q$}_S)/\partial\mbox{\boldmath$\sf q$}_S$ is detailed in \ref{sec:kinematic_linearization_1}. 

\begin{remark}
It should be noted that when the track structure is included in the finite element model, the kinematic constraint of the ``moving node'' $m$ is established with the finite element nodes of the track subsystem, rather than those of the bridge. The corresponding formulation remains unchanged.
\end{remark}


\subsection{Constraint between the moving node and the rails} \label{sec:constraint2}
\begin{figure}[!ht]
	\centering
    \subfigure[Kinematic relationship of virtual node $m$ with the virtual nodes $r_1$ and $r_2$]{\label{fig:rails}\includegraphics[width=0.45\textwidth]{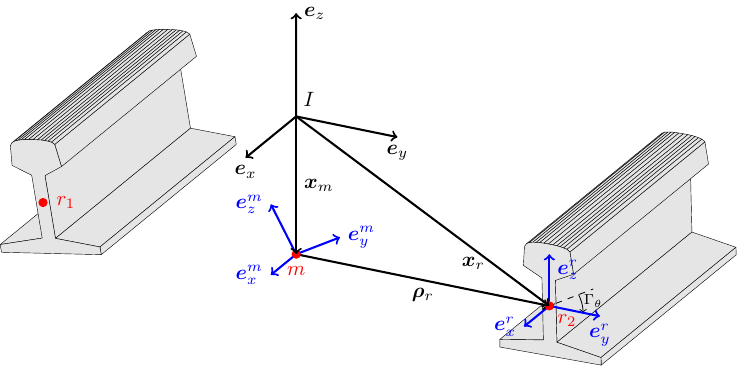}}
	\subfigure[Geometric definition of rail surface]{\label{fig:rail-extrusion}\includegraphics[width=0.45\textwidth]{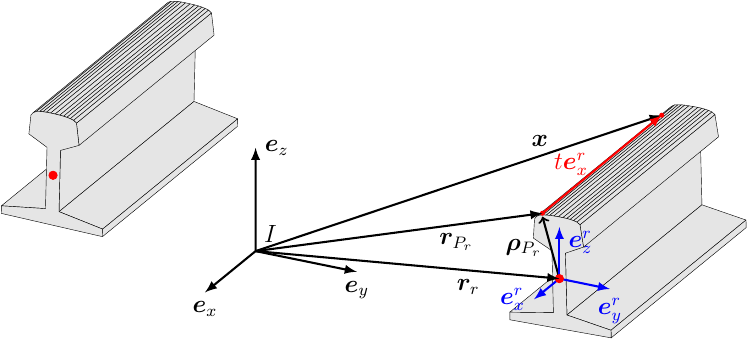}}

	\caption{Sketch of the rail definition}
\end{figure}

As previously introduced at the beginning of this section, two additional restrictions are established to define the kinematic relationship between the virtual node $m$ and the virtual nodes $r_1$ and $r_2$ shown in Fig.~\ref{fig:rails}. Since the constraint equations between $m$ and $r_1$ and the one between $m$ and $r_2$ are equivalent, they will henceforth be generalized and formulated as a constraint between virtual nodes $m$ and $r$.

\par 
Analogously to Eq.~\eqref{ch3:eq:excen}, the initial offsets of the virtual node $r$ with respect to $m$ are expressed as (see Fig.~\ref{fig:rails})
\begin{subequations} \label{eq:offsets_rails}
\begin{align}
 c^r_y & = \left(\bm{X}_r-\bm{X}_m\right)\cdot\bm{E}^m_y\,,\\
 c^r_z & = \left(\bm{X}_r-\bm{X}_m\right)\cdot\bm{E}^m_z\,,
\end{align}
\end{subequations}
being $\bm{X}_r$ the position vector of the rail centroid in the reference configuration, and $\bm{E}^m_y$ and $\bm{E}^m_z$ two unit vectors that orient the virtual node $m$. Consequently, the track section also refers to the undeformed configuration.

The kinematic relationship between $m$ and $r$ is (see Fig.~\ref{fig:rails}):
\begin{align}\label{eq:rest_disp4}
\bm{x}_r = \bm{x}_m + \bm{\rho}_r\,,
\end{align}
where $\bm{\rho}_r$ is the vector between $r$ and $m$.
\par 
Since the undeformed position of the track may differ from its ideal geometry,
the potential presence of irregularities in the track profile is considered.
These track irregularities can significantly affect the dynamic behavior of
railway vehicles \cite{Zhai2019}. In this work, the rail vector $\bm{\rho}_r$, which refers to both nodes $r_1$ and $r_2$, will take into account the track irregularity and is expressed as:
\begin{align}\label{eq:rest_disp5bis}
\bm{\rho}_r = \left(\Gamma_{y}(s) + c^r_y\right)\,\bm{e}^m_y
+ \left(\Gamma_{z}(s) + c^r_z\right)\,\bm{e}^m_z\,,
\end{align}
being $\Gamma_{y}$ and $\Gamma_{z}$, the track irregularities at section $s$\footnote{
, $\Gamma_{y}(s)$, and $\Gamma_{z}(s)$, as well as $\Gamma_{\theta}(s)$ referred to below, can be generated either by means of spectral density functions (as done in \cite{Claus1997}) or by real measurements of track deviations (see, for example, \cite{Antolin2013}) for both rails.}. The vectors $\bm{e}^m_\alpha$, with $\alpha=\lbrace x,\,y,\,z\rbrace$, correspond to the orthonormal basis linked to the node $m$ that describes the track orientation in the current configuration.
As in the case of the previous constraint, vectors $\bm{e}^m_\alpha$ are computed as
a rotation of $\bm{E}^m_\alpha$ (Eq.~\eqref{ch3:eq:rest_rot}), their
counterparts in the undeformed configuration.
\par 
Analogously to the previous section, the restriction Eq.~\eqref{eq:rest_disp4} can be rewritten as:
\begin{align}\label{eq:rest_disp4_bis}
 \bm{x}_r = \bm{X}_r + \bm{u}_r = \bm{X}_m + \bm{u}_m + \bm{\rho}_r\,,
\end{align}
and therefore $\bm{u}_r$ can be expressed as:
\begin{align}\label{eq:rest_disp5}
  \bm{u}_r = \bm{X}_m + \bm{u}_m + \bm{\rho}_r - \bm{\bm{X}}_r\,.
\end{align}
\par 
For the rotational constraint between $m$ and $r$, it is necessary to consider the possible rotation $\Gamma_\theta(s)$ of the rail around the longitudinal axis $\bm{e}^r_x$. This rotation could be induced by the track irregularities or could be due to the track cant. Thus, the rotation tensor of the rail $\bm{\Lambda}_r=\bm{\Lambda}(\bm{\phi}_r)$ is computed as the composition of two rotations through
\begin{align}\label{eq:rest_rot3}
\bm{\Lambda}_r = \bm{\Lambda}_m\,\bm{\Lambda}_\Gamma\,,
\end{align}
being $\bm{\Lambda}_\Gamma=\bm{\Lambda}(\Gamma_\theta \bm{e}^b_x)$. Therefore, the rail rotation vector $\bm{\phi}_r$ can be extracted from the tensor $\bm{\Lambda}_r$ by means of Spurrier's algorithm \cite{Spurrier1978}. 
\begin{align}\label{eq:rest_rot4}
\bm{\phi}_r = \text{spurrier}\left[\bm{\Lambda}_m\,\bm{\Lambda}_\Gamma\right]\,.
\end{align}
\par 
As in the previous case, restriction Eqs.~\eqref{eq:rest_disp5} and \eqref{eq:rest_rot4} can be formulated analogously to Eq.~\eqref{eq:restricciongeneral}, with the coordinates of node $r$ 
being the dependent ones and the coordinates of node $m$ being the independent ones.
\begin{align}\label{eq:restriccion2}
\mbox{\boldmath$\sf \Phi$}^r_{SV}= \mbox{\boldmath$\sf q$}_r
- \mbox{\boldmath$\sf\Xi$}^r(\mbox{\boldmath$\sf q$}_m)
= \mbox{\boldmath$\sf0$}\,.
\end{align}
As before, the linearization of the restriction is
\begin{align}\label{eq:restriccion3}
\delta\mbox{\boldmath$\sf q$}_r &=
\frac{\partial\mbox{\boldmath$\sf\Xi$}^r(\mbox{\boldmath$\sf
q$}_m)}{\partial\mbox{\boldmath$\sf q$}_m}\,\delta\mbox{\boldmath$\sf
q$}_m\,,
\end{align}
that is explained in detail in the \ref{sec:kinematic_linearization_2}.
\par 
It is worth noting that for this work, the general-purpose nonlinear finite element software Abaqus \citep{Abaqus610} has been used. Abaqus is a powerful tool for discretizing structures, offering a wide variety of element and material formulations, as well as multibody capabilities that enable the modeling of simplified vehicles. The kinematic constraints described in this section have been implemented in Abaqus through user subroutines for MultiPoint Constraints (MPC). The implementation of these routines is discussed in \ref{sec:mpc}.
Nevertheless, while the implementation used in this work is based on Abaqus, the proposed methodology is agnostic to the finite element software used, and similar implementations would be possible in other finite element software that supports user subroutines, such as Ansys.
\section{Wheel-rail contact model} \label{sec:contact}
The wheel-rail contact phenomenon involves both rolling and sliding interactions between two three-dimensional bodies with complexly profiled surfaces. As a result, accurately identifying the contact points and solving the associated wheel-rail interaction problem is highly challenging due to the presence of strong nonlinearities. To address this complexity, a structured approach consisting of three sequential and decoupled steps is proposed: i) determination of the contact points between wheel and rail surfaces (\ref{subsect:4.1}); ii) determination of the contact area and the normal contact forces (\ref{subsect:4.2}); and iii) evaluation of the creepages and determination of the distribution of the tangential creep forces and moments distribution (\ref{subsect:4.3}). A systematic methodology is developed to compute the contact forces between the wheel and the rail, as outlined in the following subsections. 

\subsection{Detection of contact zones between wheel and rail surfaces}
\label{subsect:4.1}
The prediction of the contact points between the wheels and rail surfaces is determined online for each time step during the dynamic analysis.
\subsubsection{Intersection between the wheel and the rail}
\begin{figure}[!ht]
	\centering

    \subfigure[Wheel surface definition]{\label{fig:wheel-surf}\includegraphics[width=0.4\textwidth]{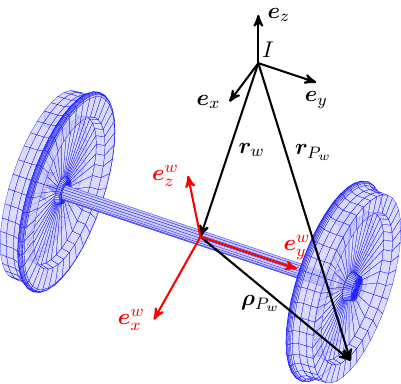}}
    \subfigure[Intersection between the straight line and truncated cone]{\label{fig:intersection}\includegraphics[width=0.55\textwidth]{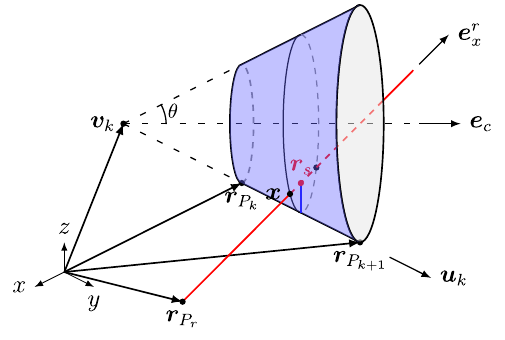}}
    \subfigure[Determination of the maximum penetration $h$ and the normal vector $\bm{n}_k$ to the cone surface]{\label{fig:penetration}\includegraphics[width=0.4\textwidth]{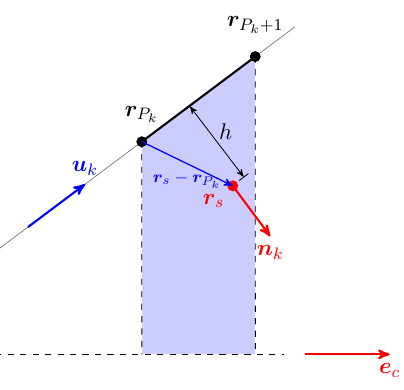}}    
	\caption{Intersection between the wheel and the rail}
\end{figure}

The wheel and rail bodies are modeled in a fully three-dimensional form, reconstructed from their actual plane profiles. The rail surface is generated by extruding its profile linearly along the longitudinal direction of the track (see Fig.~\ref{fig:rail-extrusion}), under the assumption that the rail is perfectly straight and free from surface imperfections. Since the rail profile consists of a discrete set of points connected by straight line segments, the resulting rail surface can be represented as a series of parallel lines. The discrete profile points define these lines and extend along the rail's longitudinal direction. The wheel is assumed to be perfectly smooth and free from any damage or imperfections; however, it would be straightforward to incorporate any modifications in the longitudinal or transverse directions. The surface is generated by revolving its profile around the wheel's central axis. Since the wheel profile consists of a discrete set of points connected by straight line segments, the resulting wheel surface can be interpreted as the envelope of truncated conical surfaces formed between consecutive points, as illustrated in Fig.~\ref{fig:wheel-surf}. 

To locate the contact points between the wheel and rail surfaces, it is necessary to determine the intersection of the rail's extruded lines with the truncated cones of the wheel. It is assumed that a straight line from the extrusion of the rail profile intersects the truncated cone of a wheel at a point $\bm{x}$, as illustrated in Fig.~\ref{fig:intersection}. For this to occur, the following conditions must be satisfied:
\begin{align}
	\cos \theta = \dfrac{\bm{x} - \bm{v}_k}{||\bm{x} - \bm{v}_k ||} \cdot \bm{e}_c \label{eq:cos1} \\
	[(\bm{x} - \bm{r}_{P_k})\cdot \bm{e}_c][(\bm{x} - \bm{r}_{P_{k+1}})\cdot \bm{e}_c] \leq 0 \label{eq:cone}
\end{align}
where $\theta$ is the half-angle of the truncated cone (see Fig.~\ref{fig:intersection}).

The first condition \eqref{eq:cos1} is to impose the point $\bm{x}$ on the surface of the cone, and the second condition \eqref{eq:cone} is to restrict the point $\bm{x}$, which must be on the limit of the truncated cone formed by the two points $\bm{r}_{P_k}$ and $\bm{r}_{P_{k+1}}$. The first condition \eqref{eq:cos1} can be reformulated by squaring both of its sides:
\begin{align}
	\cos^2 \theta ||\bm{x} - \bm{v}_k ||^2 = [(\bm{x} - \bm{v}_k )\cdot \bm{e}_c]^2 \nonumber \\
	\text{or} \quad (\bm{x} - \bm{v}_k )^T \bm{M} (\bm{x} - \bm{v}_k ) = 0 
	\label{eq:cos2}
\end{align}
where $\bm{M} = \bm{e}_c \otimes \bm{e}_c - \cos^2 \theta \bm{1}$ and $\bm{1}$ are the unit tensors. Replacing $\bm{x}=\bm{r}_{P_r} + t \bm{e}_x^r$ in the Eq.~\eqref{eq:cos2}, one obtains a quadratic equation:
\begin{align}
	at^2  + bt +c = 0
	\label{eq:t}
\end{align}
where $a = \bm{e}_x^r \cdot \bm{M} \bm{e}_x^r$, $b = \bm{e}_x^r \cdot \bm{M} (\bm{r}_{P_r} - \bm{v}_k)$, and $c=(\bm{r}_{P_r} - \bm{v}_k) \cdot \bm{M} (\bm{r}_{P_r} - \bm{v}_k)$, being the two possible solutions are $t_{+/-} = \left(-b \pm \sqrt{a^2-4bc}\right) / 2a$.


\noindent If $a = 0$, the rail line is parallel to the generatrix of the truncated cone, and there is a unique intersection point between both. In contrast, $a \neq 0$, three different cases can exist: i) the rail line does not intersect with the truncated cone ($a^2-4bc <0 $); ii) the rail line is tangent to the surface of the truncated cone  ($a^2-4bc = 0 $); and iii) the rail line intersects the surface of the truncated cone at two points ($a^2-4bc > 0 $).

Once the potential intersections between the straight line and the truncated cone are identified, and if there are two distinct points of intersection, it is necessary to determine the maximum penetration of this line. To calculate this value, an assumption is made that the maximum penetration corresponds to the distance from the midpoint between the two intersection points to the surface of the truncated cone (see Fig.~\ref{fig:penetration}). This assumption is exact when there is no yawing rotation of the wheelset. Additionally, it serves as a very accurate approximation, since the relative rotation between the track axis and the wheelset axis is typically small. The midpoint can be expressed as:
\begin{align}
	\bm{r}_s = \bm{r}_{P_r} - \dfrac{b}{2a} \bm{e}_x^r
\end{align}
Consequently, the maximum penetration can be obtained as follows (see Fig.~\ref{fig:penetration}):
\begin{align}
	h = ||\bm{u}_k \times (\bm{r}_s - \bm{r}_{P_k})||
\end{align}
and the normal vector $\bm{n}_k$ to the cone surface is determined as:
\begin{align}
	\bm{n}_k = \bm{u}_k \times (\bm{e}_c \times \bm{u}_k)
\end{align}
\subsubsection{Determination of multiple contact zones}

\begin{figure}[h!tp]
 \centering
 \subfigure[Two separated contact zones.]
	{\label{fig:cont1}\includegraphics[width=0.3\textwidth]{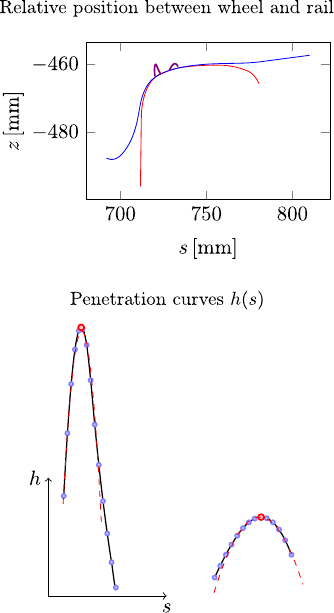}} \hfill
 \subfigure[Two joined contact zones.]
	{\label{fig:cont2}\includegraphics[width=0.3\textwidth]{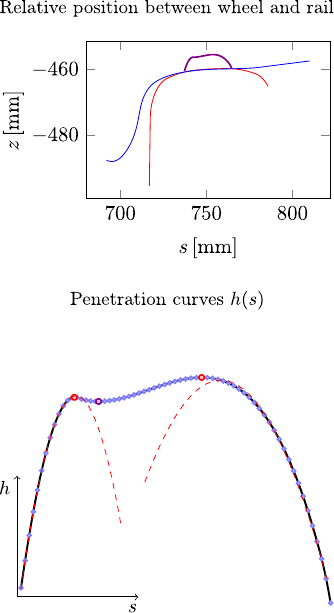}} \hfill
 \subfigure[Three joined contact zones.]
	{\label{fig:cont3}\includegraphics[width=0.3\textwidth]{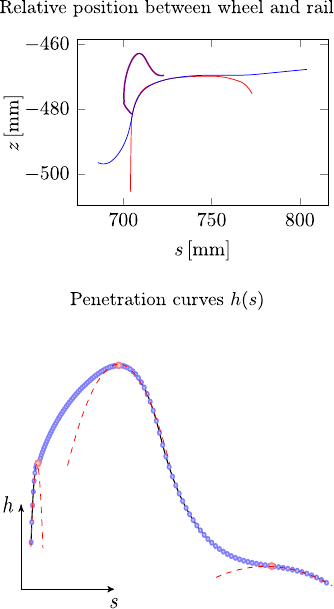}}
 \caption{Different contact zones between the wheel and rail surface. Quadratic approximation depicted in a dashed red line. Local maxima are depicted as red dots}
 
\end{figure}

For each straight line of the rail surface that intersects the wheel surface, a midpoint $\bm{r}_{si}$ and its corresponding penetration $h_i$ can be determined, allowing the identification of several contact zones (see Fig.~\ref{fig:cont1}, in which two contact zones are shown). To formulate the normal contact forces, it is essential to identify the maximum penetration and its associated position within each contact zone. A straightforward method is to find the maximum absolute value from the set of penetrations $\{h_i\}$. However, this approach introduces significant issues: as the wheel moves along the rail, even small shifts in the sampling grid can produce sudden jumps in the identified contact-point location, abrupt variations in the computed normal and tangential forces, and oscillatory behaviour in the predicted contact patch. This effect, commonly known as grid locking or discretization noise, arises from relying solely on discrete penetration values. To mitigate these problems, cubic-spline interpolation is used to determine the maximum penetration from the computed set $\{h_i\}$. The spline provides a smooth, continuously differentiable curve that passes through the sampled penetrations, allowing the maximum to be located by analysing the continuous interpolant rather than the discrete samples. This enables the algorithm to estimate the exact position of the peak within each interval and to compute a more accurate maximum penetration. To do so, a parameter $s$ is introduced, ensuring that each point $\bm{r}_{si}$ corresponds to a specific value of $s_i=\bm{r}_{si} \cdot \bm{e}_y^w$.

In the conventional contact formulation based on Hertz theory \cite{Hertz1882}, it is assumed that the surfaces in contact are described by quadratic functions of two variables at the interface. This assumption implies that the penetration curve follows a second-order polynomial form (see Fig.~\ref{fig:cont1}). However, this is not always the case. In some situations, the penetration curve deviates from a parabolic shape, as illustrated in Fig.~\ref{fig:cont2}. This behavior typically occurs when the wheelset is centered on the track \cite{Pascal1993}, and two distinct closed contact zones blend. In these cases, the number of contact zones is determined by counting the maximum penetrations of the local area (red dots in Fig.~\ref{fig:penetration} and Fig.~\ref{fig:cont2}).

It is also important to note that the local minima and maxima may occur very close together, potentially merging with the inflection point between them. This proximity can lead to convergence difficulties when solving the equations of the dynamically coupled system. As a result, special care must be taken to accurately identify all local minima, maxima, and inflection points.   

\subsection{Normal wheel-rail contact modeling}
\label{subsect:4.2}
As mentioned above, when the wheel and rail are in contact, multiple contact zones can form at the interface, which can be either separated (see Fig.~\ref{fig:cont1}) or joined (see Fig.~\ref{fig:cont2}). Consequently, approaches based on a single elliptical contact model using Hertz' theory \cite{Hertz1882} are not appropriate for these scenarios. To address all variations in wheel-rail contact and achieve accurate results, this study uses the multi-Hertzian contact method \cite{Pascal1993, Piotrowski2005}. This method assumes that each contact zone meets the necessary conditions of Hertz' theory \cite{Hertz1882, Johnson1987}, so that solving the normal contact problem involves determining a set of ellipses. The virtual compression in each contact zone $\delta_i$ is first approximated using methods established in previous studies \cite{Ayasse1987, Kik1996, Linder1997, Piotrowski2005}.
\begin{align}
\delta_i = h_i  - 0.45 h_0
\end{align}
where $h_i$ is the corresponding maximum penetration of the $i$-th contact zone, $h_0$ is the maximum penetration of the contact area (the absolute value, see Fig.~\ref{fig:cont3}) that corresponds to the ellipse 0.
Once the virtual compression is calculated, the normal contact force, the semi-axes, and the elliptical pressure distribution at each contact zone are determined considering Hertzian theory locally \cite{Hertz1882, Johnson1987, Shabana2008}.

\begin{figure}[h!pt]
	\centering
	\includegraphics[width=0.5\textwidth]{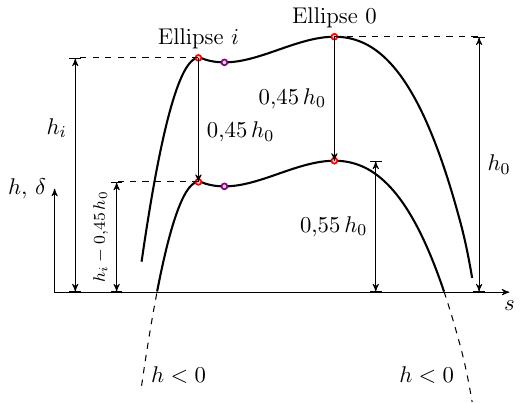}
	\caption{Approximation of the virtual compression at the contact zone, adapted from \cite{Ayasse1987, Piotrowski2005}}
	\label{fig:cont3}
\end{figure}

\subsection{Tangential wheel-rail contact modeling}
\label{subsect:4.3}
In this study, the Kalker's USETAB \cite{Kalker1996} is employed to predict the wheel-rail contact creep forces. During dynamic analysis, the longitudinal component $F_x$ and the lateral component $F_y$ of the creep force, together with the spin creep moment $M_z$, are obtained by linearly interpolating the values of the precalculated creep force stored in the USETAB lookup multi-dimensional table. These precalculated values were generated beforehand using the CONTACT software \cite{Vollebregt2012}, which is based on the \textit{Exact Three-Dimensional Rolling Contact Theory} \cite{Kalker1986,Kalker1986a}. The input provided to USETAB includes the reduced creepages and the ratio of the semi-axes $a/b$ of the contact ellipse. 

The reduced creepages are defined as follows:
\begin{subequations}
\begin{align}
	\tilde{\zeta}_{x} &= \dfrac{abG}{3 \mu N} c_{11} \zeta_x \\
	\tilde{\zeta}_y &= \dfrac{abG}{3 \mu N} c_{22} \zeta_y \\
	\tilde{\varphi} &= \dfrac{(ab)^{3/2}G}{\mu N} c_{23} \varphi
\end{align}
\end{subequations}
where $\zeta_x,\ \zeta_y,\ \text{and} \ \varphi $ are the longitudinal, lateral, and spin creepages in the global coordinate system, respectively. $G$ is the shear modulus of the material of the wheel and rail, $\mu$ is the friction coefficient between the wheel and the rail. $c_{11},\ c_{22},\ \text{and}\ c_{23}$ are the Kalker's coefficients \cite{Kalker1990} that depend on the ratio $a/b$ and the Poisson coefficient $\nu$. The longitudinal, lateral, and spin creepages are given by:
\begin{subequations}
	\begin{align}
		\zeta_x &= \dfrac{\bm{v}_w - \bm{v}_r}{v} \cdot \bm{e}_x \\
		\zeta_y &= \dfrac{\bm{v}_w - \bm{v}_r}{v} \cdot \bm{e}_y \\
		\varphi &= \dfrac{\bm{w}_w - \bm{w}_r}{v} \cdot \bm{e}_z
	\end{align}
\end{subequations}
in which, $\bm{v}_w,\ \bm{w}_w$ and $\bm{v}_r,\ \bm{w}_r$ are the vectors of the translational and angular velocities of the wheel and the rail at the contact point, respectively, $\bm{e}_x,\ \bm{e}_y,\ \text{and}\ \bm{e}_z$ are the unit base vectors of the contact point plane, and $v$ is the magnitude of the wheel velocity along the longitudinal direction, which is the same value as the train velocity.

\section{Numerical validation and applications} \label{sec:examples}

To highlight the main advantage of the proposed method—namely its unified finite element formulation in which wheel–rail contact kinematics and contact forces are derived and implemented in a fully consistent FEM framework—two numerical applications are presented.
\begin{itemize}
	\item The first application is the well-known Manchester Contact Benchmark reported by \cite{Shackleton2006, Shackleton2008}, which validates the accuracy of the proposed contact model against established industry standards.
	\item The second application investigates the dynamic interaction between a single railway coach and a continuous girder bridge under wind loading, accounting for both deck torsion and bending. This case study is intended to demonstrate the full capabilities, numerical stability, and practical applicability of the proposed unified FEM formulation for complex coupled vehicle–track–bridge systems, rather than to provide an exhaustive benchmark against all existing wheel–rail contact models.
\end{itemize}
Together, these examples demonstrate how the proposed approach provides a consistent and efficient implementation of wheel–rail interaction within a single unified computational framework, eliminating the need for additional interface assumptions and enabling fully integrated vehicle–track–bridge dynamic analyses. While a comprehensive quantitative comparison of computational efficiency and accuracy with existing models is beyond the scope of the present study and will be addressed in future work, the results clearly highlight the methodological shift toward a more integrated, consistent, and scalable framework for large-scale railway infrastructure dynamics.

\subsection{Manchester Contact Benchmarks}

\begin{figure}[h!tp]
\centering
\includegraphics[width=0.6\textwidth]{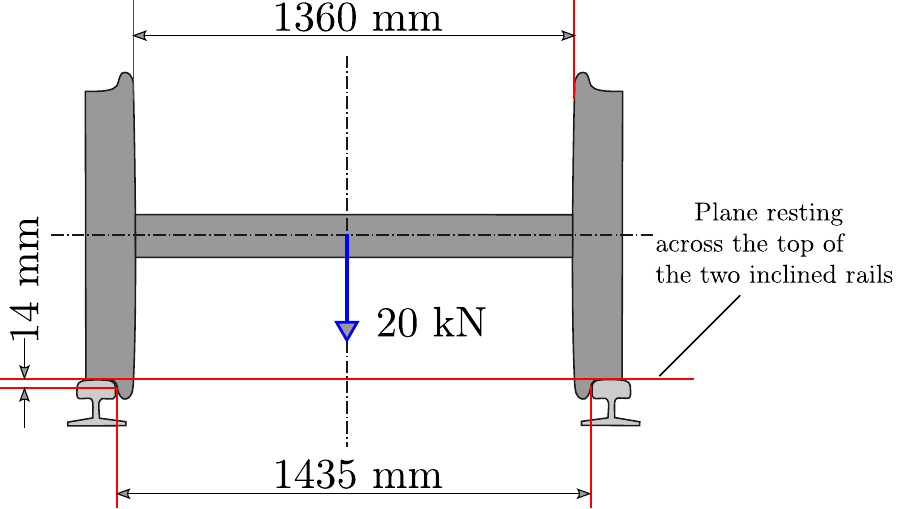}
\caption{Main geometric properties of the wheel and rail (adapted from \cite{Shackleton2008})}
\label{fig:MCB}
\end{figure}

The \textit{Manchester Contact Benchmark} proposed by \cite{Shackleton2006, Shackleton2008} is used to validate the wheel-rail contact model. The benchmark simulation case, Case A \cite{Shackleton2008}, aims to compare predicted results, such as contact size, normal and tangential stress distributions, and creep forces, using different contact models.  This simulation involves a single wheelset with a prescribed lateral displacement, yaw angle, wheel load, velocity, and friction coefficient. Throughout Case A, the following assumptions are made: the real S1002 wheel and UIC60 rail profiles are used; a vertical axle load of $20~\text{kN}$ is applied in the center of the wheelset axle (see Fig.~\ref{fig:MCB}).  

To include all aspects of interest, Case A is split into four sub-cases, depending on the type of the concerned contact and the prescribed displacements. The four sub-cases are:
\begin{itemize}
\item A1-1: a static wheelset considers only normal contact with prescribed lateral displacement. 
\item A1-2: a static wheelset considers only normal contact with prescribed lateral displacement and yaw angle.
\item A2-1: quasi-static tangential contact with prescribed lateral displacement.
\item A2-2: quasi-static tangential contact with prescribed lateral displacement and yaw angle.
\end{itemize} 

The prescribed lateral displacement varies from 0 to $10~\text{mm}$ with an increment of $0.5~\text{mm}$. For the angle of yaw, its value varies from $0~\text{mrad}$ to $24~\text{mrad}$ with an increment of $1.2~\text{mrad}$. For the case that assesses the tangential contact problem, the wheelset runs on a straight track at a constant speed of 2 m/s. In each position, the wheel-rail contact will be evaluated. The results obtained with the wheel-rail contact model proposed in this study for each simulation will be compared with the results presented by \cite{Shackleton2008}. 


\subsubsection{Contact positions and rolling radius difference}
\begin{figure}[h!tp]
	\centering
    \includegraphics[width=\linewidth]{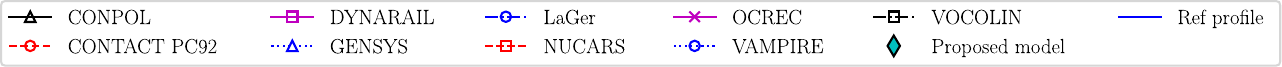}
	\subfigure[Left rail]{\includegraphics[width=0.45\textwidth]{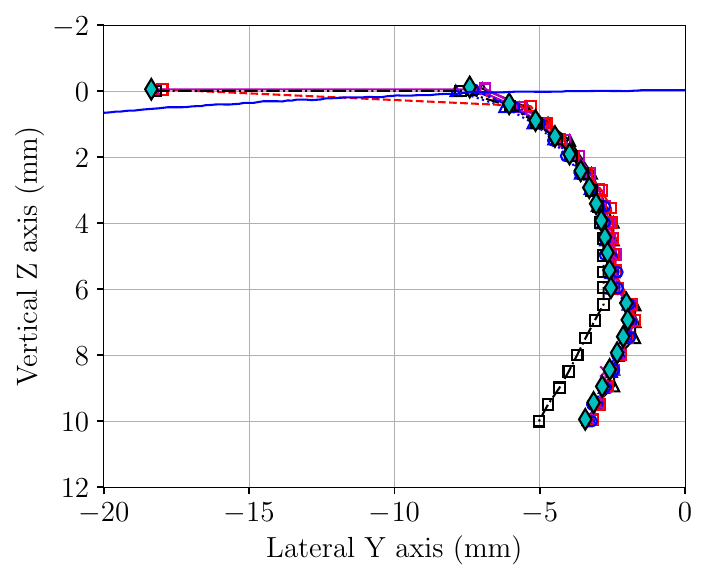}}
	\subfigure[Right rail]{\includegraphics[width=0.45\textwidth]{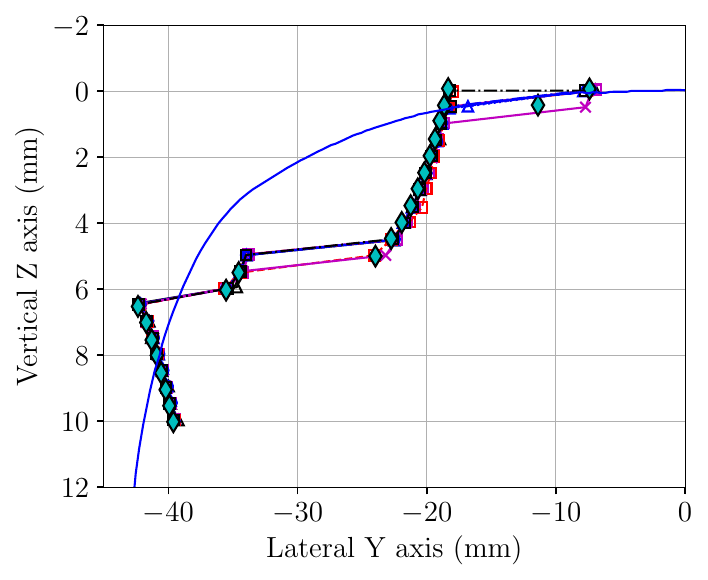}}
	\caption{Rail lateral contact positions for Case A1.1, defined in the local rail coordinate. The solid blue line with no markers shows the rail profile. (adapted from \cite{Shackleton2008})}
	\label{fig:Manchester-1.1}
\end{figure}

\begin{figure}[h!tp]
	\centering
    \includegraphics[width=\linewidth]{legend}
	\subfigure[Left wheel]{\includegraphics[width=0.45\textwidth]{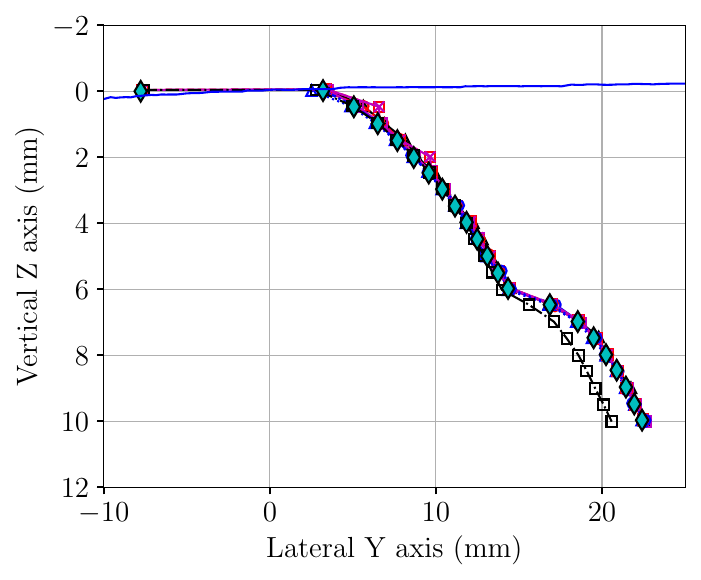}}
	\subfigure[Right wheel]{\includegraphics[width=0.45\textwidth]{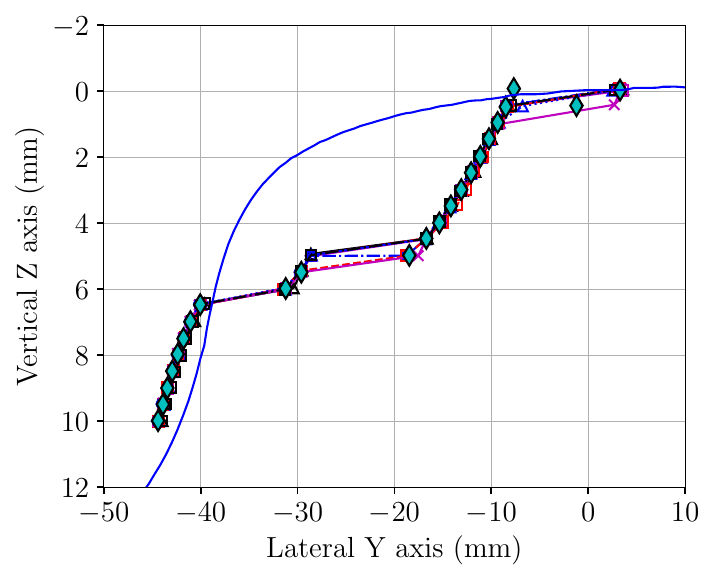}}
	\caption{Wheel lateral contact positions for Case A1.1, defined in the local rail coordinate. The solid blue line with no markers shows the wheel profile. (adapted from \cite{Shackleton2008})}
	\label{fig:Manchester-1.2}
\end{figure}

Figs.~\ref{fig:Manchester-1.1} and \ref{fig:Manchester-1.2} show the left and right contact positions defined in the local rail and wheel coordinate systems for the sub-case A1-1. It can be seen that the results obtained in this study are similar to those obtained by different contact models reported in \cite{Shackleton2008}. 


\begin{figure}[h!bt]
	\centering
    \includegraphics[width=\linewidth]{legend}
	\includegraphics[scale=0.6]{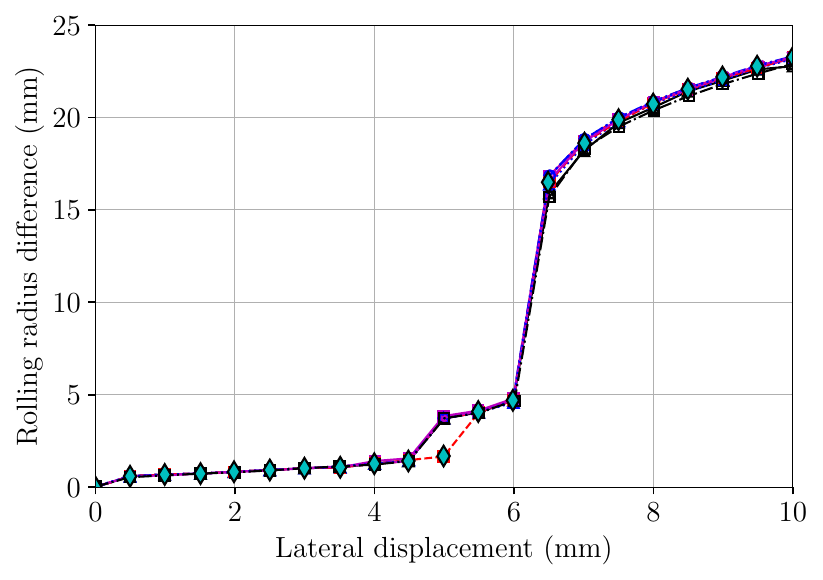}
	\caption{Difference in rolling radius for sub-case A1-1 (adapted from \cite{Shackleton2008})}
	\label{fig:Manchester-2}
\end{figure}

The difference in rolling radius between the two wheels for subcase A1-1 is shown in Fig.~\ref{fig:Manchester-2}. It can be seen again that there is consistency between the predicted results of this study and those obtained by \cite{Shackleton2008}. Furthermore, it can be observed that when the flange contact has occurred (lateral displacement greater than $6~\text{mm}$), the rolling radius difference increases sharply, leading to a strong variation of longitudinal creepages, as shown in Fig.~\ref{fig:Manchester-2}.  

\subsubsection{Contact angles}
\begin{figure}[h!bt]
	\centering
    \includegraphics[width=\linewidth]{legend}
	\subfigure[Left contact]{\includegraphics[width=0.45\textwidth]{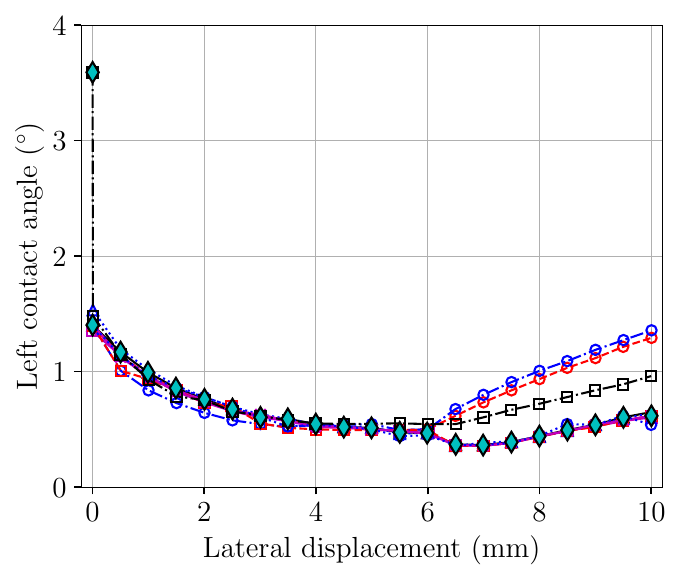}}
	\subfigure[Right contact]{\includegraphics[width=0.45\textwidth]{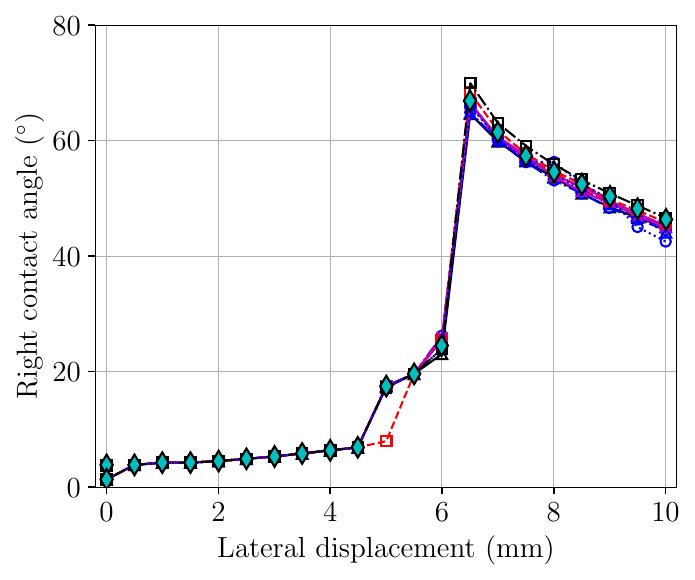}}
	\caption{Contact angle for wheel-rail interfaces for sub-case A1-2 (adapted from \cite{Shackleton2008})}
	\label{fig:Manchester-3}
\end{figure}

Fig.~\ref{fig:Manchester-3} shows the contact angle in the left and right wheel-rail interfaces for subcase A1-2. It can be seen that the results obtained in this study are closely aligned with those of CONPOL, GENSYS, and VAMPIRE. The proposed model effectively captures the sharp variation in the contact angle on the right wheel caused by flange contact when the lateral displacement exceeds $6~\text{mm}$. The maximum value is reached for a lateral displacement of $6.5~\text{mm}$.

\subsubsection{Longitudinal creepages}
\begin{figure}[!ht]
	\centering
    \includegraphics[width=\linewidth]{legend}
	\subfigure[Left contact]{\includegraphics[width=0.45\textwidth]{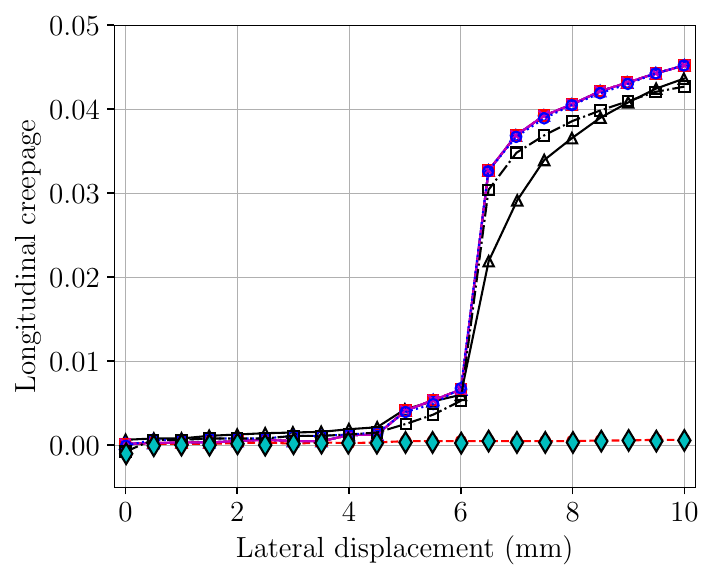}}
	\subfigure[Right contact]{\includegraphics[width=0.45\textwidth]{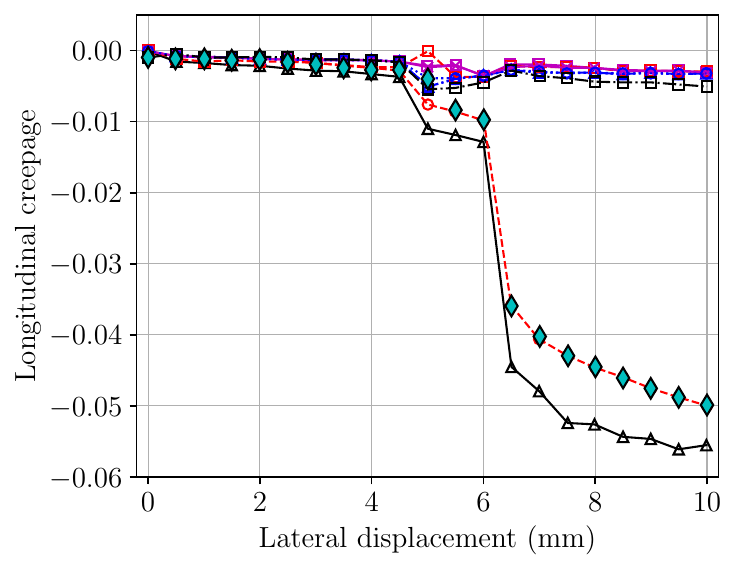}}
	\caption{Longitudinal creepages at the left and right contacts for sub-case A2-1 (adapted from \cite{Shackleton2008})}
	\label{fig:Manchester-5}
\end{figure}

The longitudinal creepages at the left and right contacts for sub-case A2-1 are depicted in Fig.~\ref{fig:Manchester-5}. As noted in \cite{Shackleton2008}, there is a significant variation in the longitudinal creepages predicted by different codes. The results obtained in this study align with those obtained by CONTACT PC92 model.
The difference between the longitudinal creepages on the left and right is consistent across all models.
\subsubsection{Spin creepages}
\begin{figure}[!ht]
	\centering
    \includegraphics[width=\linewidth]{legend}
	\subfigure[Left contact]{\includegraphics[width=0.45\textwidth]{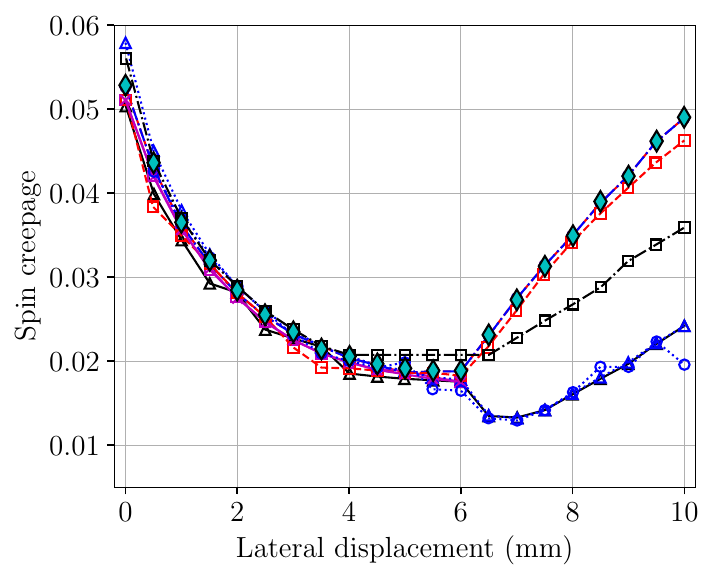}}
	\subfigure[Right contact]{\includegraphics[width=0.45\textwidth]{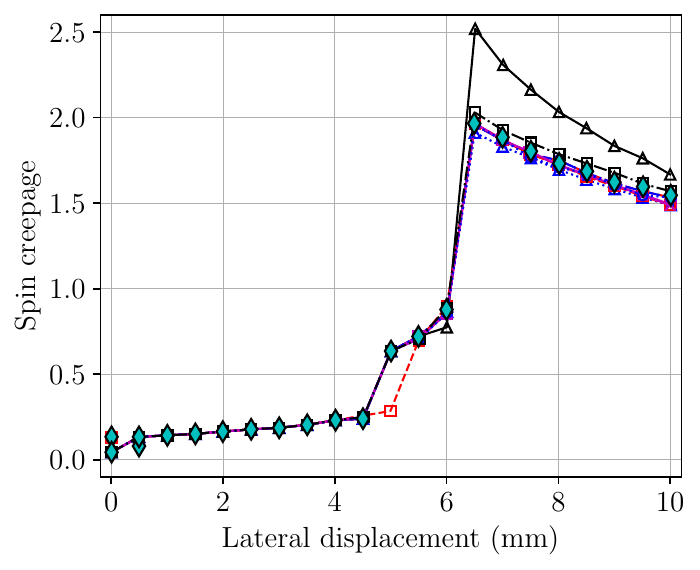}}
	\caption{Spin creepages at the left and right contacts for sub-case A2-2 (adapted from \cite{Shackleton2008})}
	\label{fig:Manchester-4}
\end{figure}

The spin creepages in this study at the left and right contact interfaces for subcase A2-2 are compared to others, as shown in Fig.~\ref{fig:Manchester-4}. As the spin creepages depend directly on the contact angle, their similarity is also observed (see Fig.~\ref{fig:Manchester-3}). 

\subsubsection{Lateral creepages}
\begin{figure}[h!bt]
	\centering
    \includegraphics[width=\linewidth]{legend}
	\includegraphics[scale=0.6]{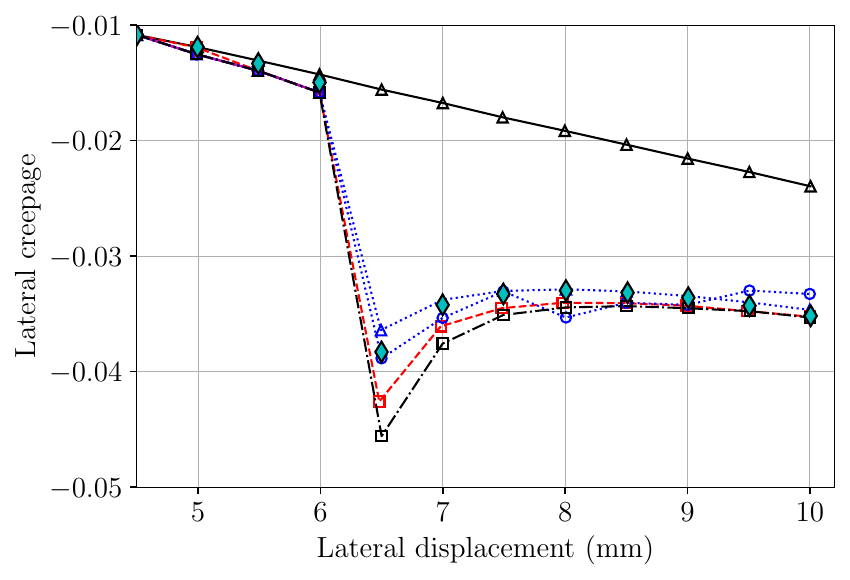}
	\caption{Lateral creepages at the right contact for sub-case A2-2, defined in the local contact patch coordinate system (adapted from \cite{Shackleton2008})}

	\label{fig:Manchester-7}
\end{figure}

In the prediction of lateral creepages, the wheel-rail model proposed in this study is in good agreement with the trend of most codes, as shown in Fig.~\ref{fig:Manchester-7}. CONPOL model is the only one that has predicted different results, as this model does not consider the effects of the angle of yaw of the wheel on the contact position, as concluded by \cite{Shackleton2008}.
\subsubsection{Contact area}
\begin{figure}[h!bt]
	\centering
    \includegraphics[width=\linewidth]{legend}
	\subfigure[Left contact]{\includegraphics[width=0.45\textwidth]{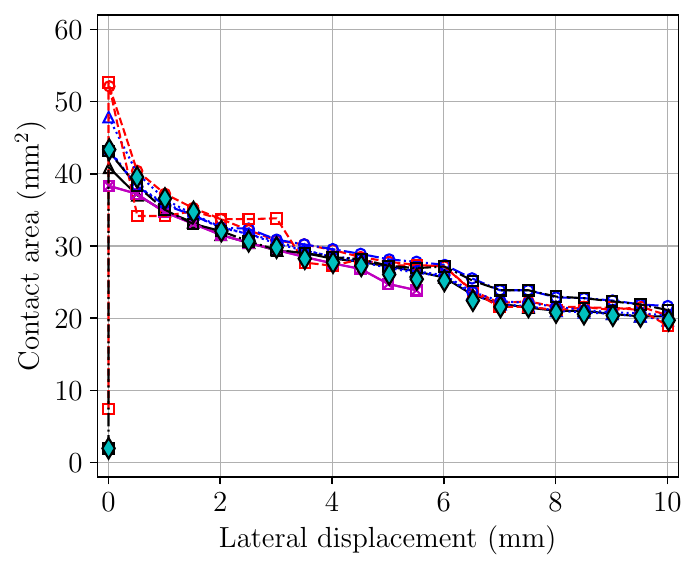}}
	\subfigure[Right contact]{\includegraphics[width=0.45\textwidth]{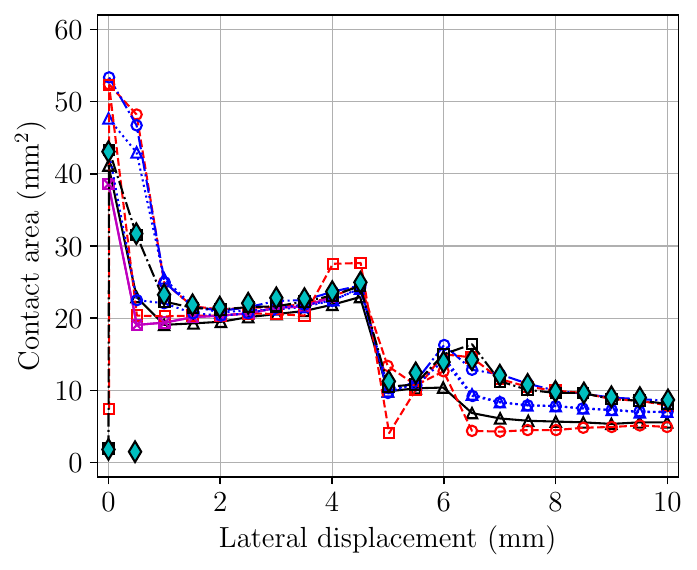}}
	\caption{Contact areas for the left and right wheel-rail contact for sub-case A1-2 (adapted from \cite{Shackleton2008}}
	\label{fig:Manchester-8}
\end{figure}

Fig.~\ref{fig:Manchester-8} shows the comparison of the contact area obtained in this study with those obtained by different contact models reported in\cite{Shackleton2008} for subcase A1-2. It can be seen that there is a clear difference in the contact area of different models after a lateral displacement of $6~\text{mm}$. The predicted contact area of this study aligns with the trend that predicts the highest values, such as LAGERS, NUCARS, and VOCOLIN models. The wheel-rail model proposed in this study, as well as those models, uses either multi-Hertzian or non-Hertzian methods for normal contact modeling and, therefore, can better represent the contact area.      
\subsection{Application for a continuous girder bridge under wind action} \label{sec:second_example}
\begin{figure}[ht]
 \centering
 \includegraphics[width=0.8\textwidth]{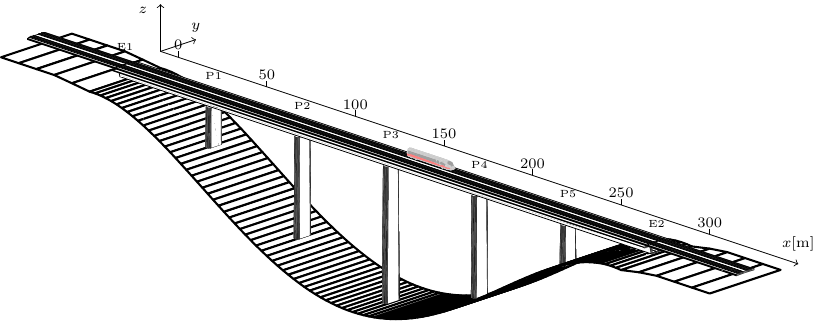}
 \caption{Scheme of the interaction example: a single vehicle crosses a
     continuous viaduct. Both the bridge and the train are subjected to a wind load that varies along
     the longitudinal coordinate $x$.}
 \label{fig:ejemplointeraccion:modelo}
\end{figure}

 In this section, as an application of the proposed model, an interaction case is studied in which a vehicle crosses a continuous girder bridge, considering both deck torsion and bending (see Fig.~\ref{fig:ejemplointeraccion:modelo}). Both the bridge and the train are subjected to aerodynamic loads that vary along the longitudinal position. These loads are applied both to the bridge and the train, deforming the structure and exciting the dynamics of the vehicle. In this example, the capabilities of the developed methodology for describing the vehicle-structure interaction are shown.
\par 
The considered structure, shown in Fig.~\ref{fig:ejemplointeraccion:modelo}, is a continuous girder bridge with a double-track section. Its total length is $300\,\text{m}$, distributed over 6 spans ($50\,\text{m}$ each). It has five piers that are $20\,\text{m}$ (P1 and P5), $40\,\text{m}$
(P2 and P4) and $60\,\text{m}$ (P3) tall. The cross section of the girder is constant along the bridge deck, and its mechanical properties are indicated in Table \ref{tab:ejemplointeraccion:props}. To demonstrate the performance of the developed models, these properties were selected to create an extremely high level of lateral bending flexibility - well beyond the serviceability limits prescribed by bridge design standards. 
\begin{table}[h!tp]
 \centering
 \begin{tabular*}{\columnwidth}{l@{\extracolsep{\fill}}cccc}
 \hline
 Part & $A\,[\text{m}^2]$ & $I_x\,[\text{m}^4]$ & $I_y\,[\text{m}^4]$ &
$J\,[\text{m}^4]$ \\
 \hline
Deck & $1.00$ & $0.95$ & $1.00$ & $2.50$ \\
P1 $\&$ P5 & $4.00$ & $2.00$ & $0.75$ & $3.00$ \\
P2 $\&$ P4 & $4.50$ & $2.50$ & $2.00$ & $4.00$ \\
P3      & $5.00$ & $3.00$ & $4.00$ & $7.00$ \\
 \hline
\end{tabular*}

 \caption{Mechanical properties of the bridge and pier sections. $A$ is the
 area of the cross-section; $I_{x}$ vertical bending inertia for the girder, and longitudinal for the piers;
 $I_{y}$ lateral bending inertia for the bridge deck, transversal for the piers; and $J$ is the torsional inertia.}
 \label{tab:ejemplointeraccion:props}
\end{table}
\par 
Both the girder and the piers have been modeled with Timoshenko beam elements (shear deformable) approximately $1\,\text{m}$ long, using a linear elastic material that corresponds to prestressed concrete in the elastic range ($E=30\,\text{GPa}$, $\nu=0.2$, and $\rho=2500\,\text{kg}/\text{m}^3$). In all the supports (abutments and deck-pier connections), the relative translations are constrained, as well as the torsional rotation. Bending rotations are allowed. In that way, when a vehicle runs along the viaduct, it produces torsional momentum in the deck section, also inducing transverse bending on the piers.
\par 
On the other hand, the vehicle considered in this example is a single car of a high-speed train that runs at a constant velocity $v=200\,\text{km}/\text{h}$ on the leeward track. The multibody model of the vehicle is composed of a car body and two bogies, each containing two wheelsets. Two suspension systems link the different
bodies: the primary suspension connects the bogies and wheelsets, while the secondary suspension links the car body to both bogies. Its main features are listed in Table \ref{tab:ejemplointeraccion:vehicle_props}.
\begin{table}[h!tp]
 \centering
 \begin{tabular*}{\columnwidth}{l@{\extracolsep{\fill}}ccl@{\extracolsep{\fill}}
 cc}
 \hline
 Item & Unit & Value & Item & Unit & Value\\
 \hline
 $m_c$ & $\text{kg}$ & $53500.0$ & $I_{c x}$ & $\text{m}^4$ & $70000.0$ \\
 $I_{c y}$ & $\text{m}^4$ & $2621500.0$ & $I_{c z}$ & $\text{m}^4$ & $2621500.0$ \\
 $m_b$ & $\text{kg}$ & $3500.0$ & $I_{b x}$ & $\text{m}^4$ & $560.0$ \\
 $I_{b y}$ & $\text{m}^4$ & $315.0$ & $I_{b z}$ & $\text{m}^4$ & $1715.0$ \\
 $m_w$ & $\text{kg}$ & $1800.0$ & $I_{w x}$ & $\text{m}^4$ & $1000.0$ \\
 $I_{w y}$ & $\text{m}^4$ & $100.0$ & $I_{w z}$ & $\text{m}^4$ & $1000.0$ \\
 $k_{x 1}$ & $\text{kN}/\text{m}$ & $120000.0$ & $k_{y 1}$ & $\text{kN}/\text{m}$ & $12500.0$ \\
 $k_{z 1}$ & $\text{kN}/\text{m}$ & $1200.0$ & $c_{x 1}$ & $\text{kN}\cdot\text{s}/\text{m}$ & $27.9$\\
 $c_{y 1}$ & $\text{kN}\cdot\text{s}/\text{m}$ & $9.0$ & $c_{z 1}$ & $\text{kN}\cdot\text{s}/\text{m}$ & $10.0$ \\
 $k_{x 2}$ & $\text{kN}/\text{m}$ & $12000.0$ & $k_{y 2}$ & $\text{kN}/\text{m}$ & $240.0$ \\
 $k_{z 2}$ & $\text{kN}/\text{m}$ & $350.0$ & $c_{x 2}$ & $\text{kN}\cdot\text{s}/\text{m}$ & $600.0$ \\
 $c_{y 2}$ & $\text{kN}\cdot\text{s}/\text{m}$ & $30.0$ & $c_{z 2}$ &
 $\text{kN}\cdot\text{s}/\text{m}$ & $20.0$ \\
 $h_c$ & $\text{m}$ & $1.4$ & $h_b$ & $\text{m}$ & $0.5$ \\
 $r_0$ & $\text{m}$ & $0.46$ & $d_c$ & $\text{m}$ & $17.375$ \\
 $d_b$ & $\text{m}$ & $2.5$ & $d_w$ & $\text{m}$ & $0.753$ \\
 \hline
\end{tabular*}

 \caption{Mechanical properties of the vehicle. The subindex $c$ corresponds to the car body, $b$ to the bogies, and $w$ to the wheelsets. Subindexes $x$, $y$, and $z$ refer to the corresponding coordinate directions. $m$ are the masses  of the bodies and $I$ are their inertias. The constants $k$ and $c$ refer to the stiffness and damping of the two suspension systems (subindex $1$ for the primary system and $2$ for the secondary). Finally, $h_c$ is the height (concerning the top of the rails) of the center of mass of the carbody; $h_b$ corresponds to the bogies; $r_0$ is the nominal rolling radius of the wheels; $d_c$ is the distance,
 along $x$, between the center of the bogies; $d_b$ is the longitudinal distance between two wheelsets of the same bogie; and $d_w$ is the semi-distance between the contact points of two wheels of the same wheelset in the rest situation.} \label{tab:ejemplointeraccion:vehicle_props}
\end{table}
\par
\begin{figure}[ht]
 \centering
 \subfigure[Chinese hat wind velocity function]
 {\label{fig:chinese_hat}
 \includegraphics{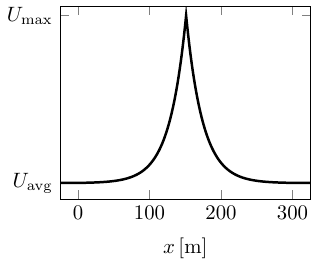}} \hspace{2cm}
 \subfigure[Sketch of the wind load components on the bridge section]
 {\label{fig:wind_loads_deck}
 \includegraphics{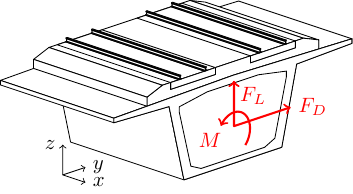}}
 \caption{Chinese hat wind velocity function and aerodynamic loads
 for a section of the bridge deck.}
\end{figure}

 
To excite vehicle and structural dynamics, the action of a wind gust was applied. Indeed, the Chinese hat wind model was deliberately adopted as a simplified but severe loading scenario. It provides a well-defined and highly transient gust that represents an extreme situation, capable of inducing large wheel–rail force variations and potential safety-critical responses. This makes it particularly suitable as a stress test for the proposed contact formulation, while avoiding the additional complexity and uncertainty introduced by stochastic wind models. It has been modeled through a Chinese hat function \citep{UNE-EN14067-6} whose average velocity is $25\,\text{m}/\text{s}$ and $35.17\,\text{m}/\text{s}$ the maximum value. The variation in wind speed corresponds to a stationary function that depends on the longitudinal position along the bridge (see Fig.~\ref{fig:chinese_hat}), with the maximum value occurring at $x=150\,\text{m}$ (pier number 3). 
Since the vehicle moves with velocity $v$, the train perceives the variation in wind velocity as a function of time. However, it is a static load on the bridge.
\par 
Neither the vertical nor the longitudinal components of wind velocity are considered; only the transverse component (along the $y$ axis according to Fig.~\ref{fig:wind_loads_deck}). Therefore, only two-dimensional aerodynamic loads appear both on the bridge and the train: they are the lift forces $F_L$ (along $z$) and the drag forces $F_D$ (along $y$), as well as the overturning moment $M$ (around $x$) (see Fig.~\ref{fig:wind_loads_deck}). These loads can be computed using the expressions \citep{Simiu1996}.
\begin{align}
    F_D &= \frac{1}{2}\,\rho_A\,A\,U^2_w\,c_D\,,\\
    F_L &= \frac{1}{2}\,\rho_A\,A\,U^2_w\,c_L\,,\\
    M &= \frac{1}{2}\,\rho_A\,A\,U^2_w\,c_M\,d\,,
\end{align}
where $\rho_A$ is the density of the air (taken as $1.25\,\text{kg}/\text{m}^3$); $A$ and $d$ are the normalized area and length of reference, respectively; $U_w$ is the wind velocity; and $c_D$, $c_L$, and $c_M$ are the corresponding aerodynamic coefficients. The expressions above apply to both the vehicle and the bridge.
\par 
The wind gust is introduced using the Chinese hat model. Its action is assumed to be instantaneous, such that the bridge deformation prior to the gust reaching its maximum intensity can be neglected. Accordingly, the aerodynamic coefficients are evaluated for a zero angle of incidence. The aerodynamic coefficients adopted for the bridge are: $c_D=0.64$, $c_L=0.34$, and $c_M=-0.31$, and they are dimensionless for $d=14\,\text{m}$ and $A=L\,d$, with $L$ being the girder section considered.
\par 
For the vehicle, $U_w$ is a composition of the wind and vehicle velocities. Thus, the wind and the train trajectory directions form a horizontal angle of attack $\beta$ (see \cite{UNE-EN14067-6} for further details).
The aerodynamic coefficients prescribed in \cite{UNE-EN14067-6} have been used for a vehicle similar to the Siemens Velaro high-speed train. The three coefficients depend on the horizontal angle of attack $\beta$.
\par 
The objective of this example is to demonstrate the capabilities of the methodology developed to couple a vehicle with a deformable structure. For this reason, the results obtained (for the \emph{flexible case}) are
compared to those computed under the assumption that the structure is completely rigid (the \emph{rigid case}). For simplicity, in this case, track irregularities have not been considered.
\par 
In all the values shown below, the value of $x$ on the abscissa axis indicates the longitudinal position of the first wheelset of the vehicle: $x=0\,\text{m}$ corresponds to the position over the first abutment, and $x=300\,\text{m}$ corresponds to the position over the second one.
\par
\begin{figure}[!ht]
 \centering
 \subfigure[Forward wheelset of the forward bogie]
 {\label{fig:ejemplointeraccion:cargas_eje1}\includegraphics{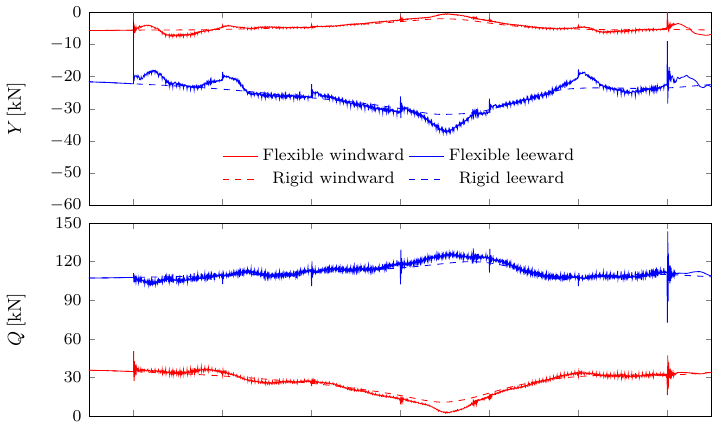}}\\
 \subfigure[Backward wheelset of the forward bogie]
 {\label{fig:ejemplointeraccion:cargas_eje2}\includegraphics{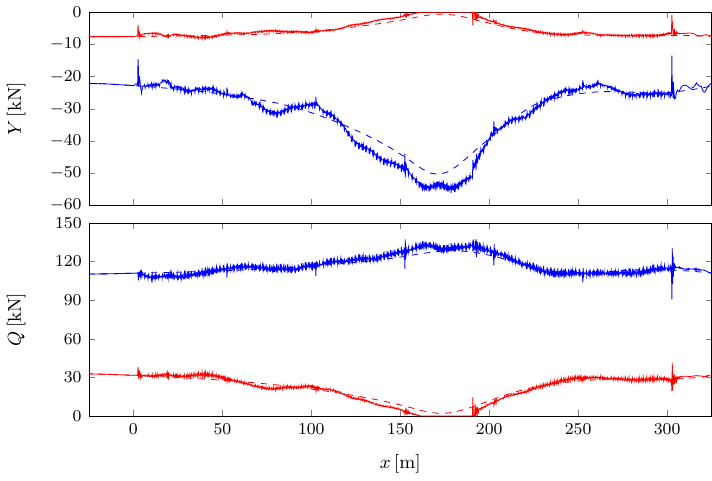}}
 \caption{Histories of the vertical $Q$ and transversal $Y$ contact loads
 for all the wheels of the forward bogie the train.}
 \label{fig:ejemplointeraccion:cargas}
\end{figure}

 
In Fig.~\ref{fig:ejemplointeraccion:cargas}, the vertical $Q$ and transverse $Y$ contact loads are shown for all the wheels of the forward bogie of the train. As can be seen in this picture, the vertical and transverse loads differ depending on whether the structure is considered flexible or rigid. In the flexible case, it can be seen that for the fourth span ($x\approx150\,\text{m}$, where the aerodynamic load is stronger and, consequently, the deformation of the bridge is also greater), the vertical loads are null; i.e., the windward wheels are not in contact with the rails. Therefore, security is more compromised in the flexible model compared to the rigid one. The differences in the force values are due not only to the flexibility of the structure but also to the fact that the track is not perfectly straight or horizontal, as the bridge is deformed under the action of the wind.
\par 

\par 
It is worth mentioning that, in the flexible case, high frequencies appear in the load histories, which are induced by the structural response. In addition, as can be seen at the entrance and exit of the viaduct (for $x=0\,\text{m}$ and $x=300\,\text{m}$), an abrupt change appears in the vertical load value. This effect is caused by the crossing of the wheelset from a rigid platform to a flexible one and vise versa. In the same way, a similar but smoother effect can be appreciated when the wheelsets cross over the pier sections (at $x=50,\,100,\,150,\,200,\,250\,\text{m}$).
\par 
In this example, the safety of the vehicle is assessed using the Nadal formula \citep{TSI2008rolling} and the offload coefficient. The Nadal formula evaluates the wheel-climb derailment risk through $\eta_N = Y/Q$. Usually, it is considered that, to guarantee safety, it must be ensured that $\eta_N<0.8$ is evaluated for each wheel of the train.
\par 
On the other hand, the offload coefficient is \citep{UNE-EN14067-6}
\begin{align}
    \eta_O = 1 - \frac{Q_{i 1}+Q_{j 1}}{2\,Q_0}\,,
\end{align}
where $Q_{i 1}$ and $Q_{j 1}$ are the vertical load values of the wheels offloaded onto the two wheelsets of the same bogie, $Q_0$ is the value of the vertical load of a vehicle wheel in a static situation. The safety limit for the coefficient, according to \cite{UNE-EN14067-6}, is $\eta_O<0.9$. Furthermore, as is well known, wheel unloading and derailment are phenomena dominated by low-frequency dynamics. Accordingly, following EN 14067-6 \cite{UNE-EN14067-6}, a fourth-order Butterworth low-pass filter with a cutoff frequency of 2 Hz is applied to the wheel–rail forces to compute the offload coefficient. For the Nadal coefficient, a different procedure specified in UIC 518 \cite{UIC518} has been adopted. This standard requires computing a moving average over a distance of 2 m from the time-domain signals of the Nadal coefficient at each wheel and axle. The resulting signal is then passed through a low-pass filter with a cutoff frequency of 20 Hz.
\begin{figure}[!ht]
 \centering
 \includegraphics{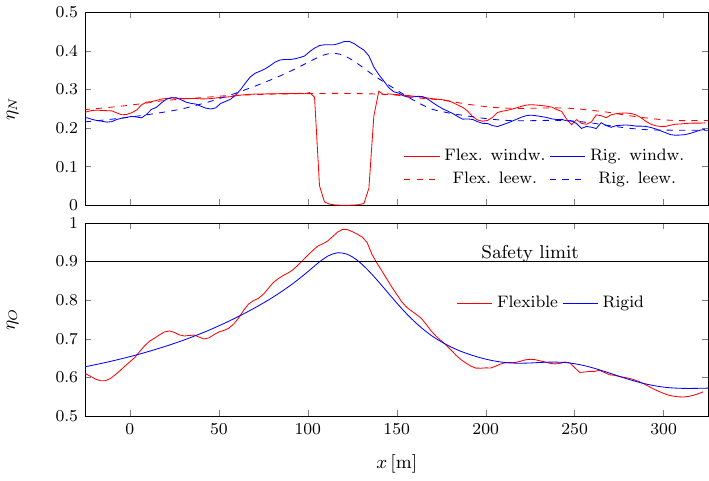}
 \caption{Histories of the safety coefficients of the vehicle relative
 to the forward bogie. The  offload coefficient $\eta_O$ corresponds to the forward
 bogie and the Nadal coefficient $\eta_N$ to the forward wheelset of the same bogie}
 \label{fig:safety_coeffs}
\end{figure}
\par 
In Fig.~\ref{fig:safety_coeffs}, the histories of $\eta_N$ and $\eta_O$ are shown for the forward bogie of the train. As can be seen, structural flexibility significantly affects the values of the safety coefficients, making the flexible case more compromising. It can be seen that the Nadal coefficient $\eta_N$ is null in section $x\approx150\,\text{m}$, because the wheel load is null $Q=Y=0$. The offload coefficient overrides the safety limit for both flexible and rigid cases.
\par
\begin{figure}[!ht]
 \centering
 \includegraphics{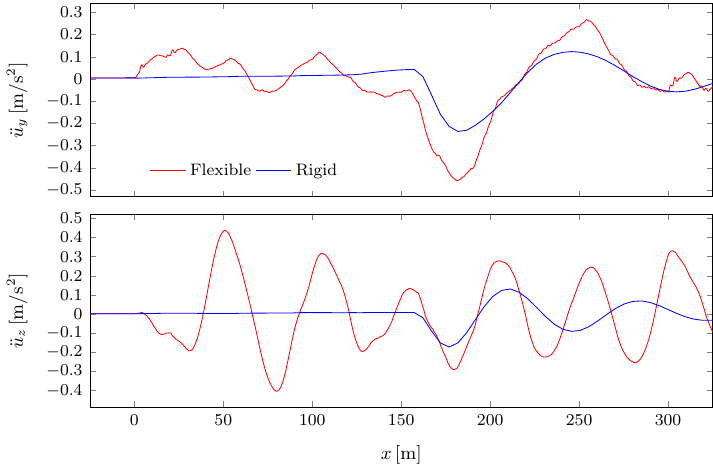}
 \caption{Histories of lateral $\ddot{u}_y$ and vertical $\ddot{u}_z$ accelerations computed on the
 center of mass of the car body of the vehicle.}
 \label{fig:ejemplointeraccion:accs_caja}
\end{figure}

 
The lateral and vertical accelerations of the center of mass of the body of the passenger car are shown in Fig.~\ref{fig:ejemplointeraccion:accs_caja}. The maximum values of the accelerations are greatly affected by the flexibility of the bridge.
\par
\begin{figure}[!ht]
 \centering
 \includegraphics{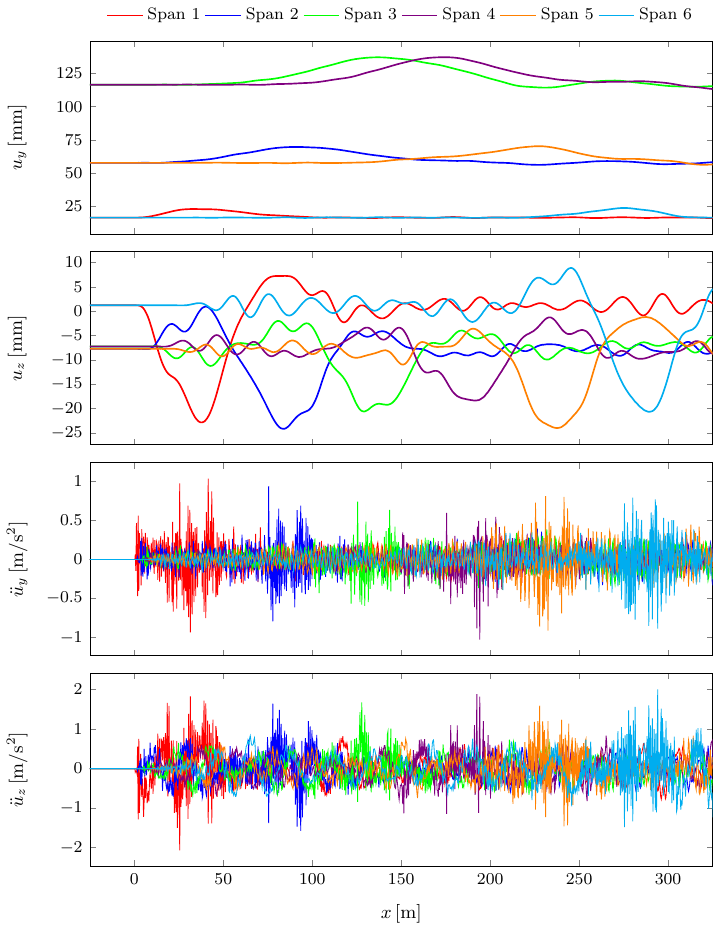}
 \caption{Histories of lateral $u_y$ and vertical $u_z$ displacements and lateral $\ddot{u}_y$ and vertical $\ddot{u}_z$ accelerations of the bridge deck computed on the center axis of the leeward track, at the top of the rails.} 
 \label{fig:ejemplointeraccion:cdv}
\end{figure}

 
Finally, in Fig.~\ref{fig:ejemplointeraccion:cdv}, vertical and lateral displacements and accelerations are shown for the center span sections of the bridge deck. Those histories are computed at the center axis of the leeward track, at the top of the rails.
\par 
The displacements are not initially null because the wind action is applied before the vehicle enters the bridge. However, the stationary nature of the wind initially causes the acceleration of the bridge deck to be null. Therefore, the accelerations are solely due to the vehicle's action.
\par 
Taking into account that a single car of a full train is considered, the bridge displacements, both lateral and vertical, are very high; the same applies to the accelerations: the bridge flexibility does not correspond to a real case.

\section{Summary and conclusions} \label{sec:conclusions}
We have proposed a methodology for establishing the kinematic relationships required for the dynamic interaction between vehicles and structures. The number of assumptions has been kept small to achieve a completely general method that facilitates every situation. The methodology is based on the use of absolute coordinates, which facilitates its integration into existing finite element and multibody dynamics software. It consists of the introduction of a set of virtual nodes that describe the track section associated with each vehicle wheelset. The position and orientation of these nodes account for the deformation of the structure and allow for introducing a richer description of the track geometry (irregularities, cant deficiency, \etc). The numerical tests performed demonstrate the accuracy of the method and allow for to envisage the capabilities and applicability of the proposed methodology.
\appendix
\section{Linearization of restrictions} \label{sec:kinematic_linearization}
In this section,  before developing the linearization of the restrictions 
\eqref{eq:restricciongeneral3} and \eqref{eq:restriccion2} 
in \ref{sec:kinematic_linearization_1} and \ref{sec:kinematic_linearization_2},
respectively, the rotation tensor and the linearization of
rotation vectors are presented in \ref{sec:rotations}.

\subsection{Rotation tensor and variation of the rotation field}
\label{sec:rotations}
According to Rodrigues formula, the rotation tensor $\bm{\Lambda}$ associated
to a rotation vector $\bm{\phi}$ can
be expressed as
\begin{align}
    \bm{\Lambda}(\bm{\phi})=\exp\left[\widehat{\bm{\phi}}\right] =
    \cos\abs{\bm{\phi}}\,\bm{1}
    + \frac{\sin\abs{\bm{\phi}}}{\abs{\bm{\phi}}}\widehat{\bm{\phi}}
    + \frac{1-\cos\abs{\bm{\phi}}}{\abs{\bm{\phi}}^2}\bm{\phi}\otimes\bm{\phi}\,,
\end{align}
where $\widehat{\bm{\phi}}$ is the skew tensor defined such as
$\widehat{\bm{\phi}}\,\bm{a}=\bm{\phi}\times\bm{a}$, $\bm{1}$ is the identity tensor
and $\otimes$ is the dyadic product of two vectors.
\par 
On the other hand, the variation of the rotation field
$\bm{\phi}$  can be written as
\begin{align}
\delta\bm{\phi} = \bm{T}(\bm{\phi})\,\delta\bm{\theta}
\end{align}
where $\bm{\theta}$ is the vector field of linearized rotations.
The operator $\bm{T}$ can be expressed as (further details can be found in
\cite{Simo1988})
\begin{align}
\bm{T}(\bm{\phi}) = 
\frac{\abs{\bm{\phi}}\,\sin\abs{\bm{\phi}}}{2\,\left(1-\cos\abs{\bm{\phi}}\right)}
\bm{1} -\frac{1}{2}\,\widehat{\bm{\phi}}+
\frac{1}{\abs{\bm{\phi}}^2}\,\left(1-
\frac{\abs{\bm{\phi}}\,\sin\abs{\bm{\phi}}}{2\,\left(1-\cos\abs{\bm{\phi}}\right)}
\right)\bm{\phi}\otimes\bm{\phi}\,.
\end{align}
And the corresponding inverse operator $\bm{H}=\bm{T}^{-1}$ allow us to write
\begin{align}
\delta\bm{\theta} = \bm{H}(\bm{\phi})\,\delta\bm{\phi}
\end{align}
being
\begin{align}
\bm{H}(\bm{\phi}) = 
\frac{\sin\abs{\bm{\phi}}}{\abs{\bm{\phi}}}\,\bm{1}
+\frac{1-\cos\abs{\bm{\phi}}}{\abs{\bm{\phi}}^2}\,\widehat{\bm{\phi}}+
\frac{1}{\abs{\bm{\phi}}^2}\,\left(1-
\frac{\sin\abs{\bm{\phi}}}{\abs{\bm{\phi}}}
\right)\bm{\phi}\otimes\bm{\phi}\,.
\end{align}

\subsection{Linearization of the constraint between the moving node and the structure} \label{sec:kinematic_linearization_1}
Below, the linearization of the constraint between the moving node $m$ and the structure, introduced
in Section \ref{sec:constraint1}, is detailed.
The variation of the rotation restriction \eqref{eq:interp_rotation} is
\begin{align}
\bm{T}(\bm{\phi}_m)\,\delta\bm{\theta}_m =
N_1(\zeta)\,\bm{T}(\bm{\phi}_1)\,\delta\bm{\theta}_1 +
N_2(\zeta)\,\bm{T}(\bm{\phi}_2)\,\delta\bm{\theta}_2\,,
\end{align}
that can be rearranged as
\begin{align}\label{eq:variation_theta}
\delta\bm{\theta}_m =
N_1(\zeta)\,\bm{H}(\bm{\phi}_m)\,\bm{T}(\bm{\phi}_1)\,\delta\bm{\theta}_1 +
N_2(\zeta)\,\bm{H}(\bm{\phi}_m)\,\bm{T}(\bm{\phi}_2)\,\delta\bm{\theta}_2\,,
\end{align}
that is the linearization of the rotation of the node
$m$ expressed as linear combination of the linearized rotations of the beam
nodes $\delta\bm{\theta}_1$ and $\delta\bm{\theta}_2$.
\par 
On the other hand, the variation of the displacement restriction \eqref{eq:rest_disp3} is
\begin{align}
 \delta\bm{u}_m =
 N_1(\zeta)\,\delta\bm{u}_1+N_2(\zeta)\,\delta\bm{u}_2 + \delta\bm{\rho}_m\,,
\end{align}
where the variations of the constant vectors $\bm{X}_m$, $\bm{X}_1$ and $\bm{X}_2$ vanish.
According to \cite{Simo1988}, the variation $\delta\bm{\rho}_m$ is:
\begin{align}\label{eq:variation_rhom}
 \delta\bm{\rho}_m = \widehat{\delta\bm{\theta}_m}\,\bm{\rho}_m = -
 \widehat{\bm{\rho}_m}\,\delta\bm{\theta}_m\,,
\end{align}
where $\widehat{\bm{\rho}_m}$, as introduced above, is the skew tensor
such as $\widehat{\bm{\rho}_m}\,\bm{a} = \bm{\rho}_m\times\bm{a}$.
Thus, considering Eqs.~\eqref{eq:variation_theta} and \eqref{eq:variation_rhom},
$\delta\bm{u}$ is expressed as linear combination of the variations
$\delta\bm{u}_1$, $\delta\bm{u}_2$, $\delta\bm{\theta}_1$ and $\delta\bm{\theta}_2$:
\begin{multline}
 \delta\bm{u}_m =
 N_1(\zeta)\,\delta\bm{u}_1+N_2(\zeta)\,\delta\bm{u}_2
 - N_1(\zeta)\,\widehat{\bm{\rho}_m}\,\bm{H}(\bm{\phi}_m)\,\bm{T}(\bm{\phi}_1)\,\delta\bm{\theta}_1\\
 - N_2(\zeta)\,\widehat{\bm{\rho}_m}\,\bm{H}(\bm{\phi}_m)\,\bm{T}(\bm{\phi}_2)\,\delta\bm{\theta}_2\,.
\end{multline}

\subsection{Linearization of the constraint between the moving node and the rails} \label{sec:kinematic_linearization_2}
In this section, the linearization of the constraint between the moving nodes $m$ and $r$
(referring indistinctly to $r_1$ and $r_2$), that was developed in Section \ref{sec:constraint2}, is presented. 
The linearization of the rotation restriction \eqref{eq:rest_rot3}
is (further details in \cite{Simo1988})
\begin{align}
\widehat{\delta\bm{\theta}_r}\bm{\Lambda}_r =
\widehat{\delta\bm{\theta}_m}\bm{\Lambda}_m\,\bm{\Lambda}_\Gamma\,,
\end{align}
where it is assumed that the rotation $\Gamma_\theta$ is constant for a
given track position $s$, thus $\delta\bm{\Lambda}_\Gamma=\bm{0}$.
Taking into account the orthogonality of the rotation tensor, the definition of
$\bm{\Lambda}_r$ (Eq.~\eqref{eq:rest_rot3}),
and right-multiplying the above expression by $\bm{\Lambda}_r^\text{T}$,
it can be written as
\begin{align}
\widehat{\delta\bm{\theta}_r} = \widehat{\delta\bm{\theta}_m}\,,
\end{align}
being the associated skew vectors
\begin{align}
\delta\bm{\theta}_r = \delta\bm{\theta}_m\,.
\end{align}
\par 
On the other hand, the variation of the displacement restriction (Eq.~\eqref{eq:rest_disp5})
is
\begin{align}
\delta\bm{u}_r = \delta\bm{u}_m + \delta\bm{\rho}_r
\end{align}
where the variations of the constant vectors $\bm{X}_m$ and $\bm{X}_r$ vanish.
Finally, the variation $\delta\bm{\rho}_r$ is computed as
\begin{align}
 \delta\bm{\rho}_r = \widehat{\delta\bm{\theta}_r}\,\bm{\rho}_r
 = - \widehat{\bm{\rho}_r}\,\delta\bm{\theta}_r
 = - \widehat{\bm{\rho}_r}\,\delta\bm{\theta}_m\,,
\end{align}
Therefore, the linearization of the displacement restriction is expressed as a linear combination of
$\delta\bm{u}_m$ and $\delta\bm{\theta}_m$:
\begin{align}
\delta\bm{u}_r = \delta\bm{u}_m - \widehat{\bm{\rho}_r}\,\delta\bm{\theta}_m\,.
\end{align}

\section{Abaqus user surboutine for constraints} \label{sec:mpc}

The finite element software Abaqus supports the possibility of extending its capabilities by allowing users 
to develop new routines for specific problems: the so-called Abaqus User Subroutines \cite{Abaqus610}. Among how the software can be extended, Abaqus offers the possibility of including new defined MultiPoint Constraints (MPC). MPC allows the establishment of kinematic constraints for the degrees of freedom (displacements and rotations) of a single node (slave node) as a function of the degrees of freedom and positions of a collection of nodes (master nodes). This Abaqus feature is perfectly suited for the implementation of the kinematic constraints proposed in this work.

Hence, the first constraint (Section \ref{sec:constraint1}) is established between the node $m$ (slave node) and all the nodes that define the structure and track (master nodes). On the other hand, for the second constraint (Section \ref{sec:constraint2}), the node $m$ is the master node, whereas the node $r$ (which refers either to $r_1$ or $r_2$) is the slave node.

In this context, two user subroutines have been developed for implementing both constraints. Their pseudocode is included in Algorithms \ref{alg:constraint_1} and \ref{alg:constraint_2}, respectively. As previously introduced, the kinematic constraints detailed in Section \ref{sec:kinematics} are set up for every wheelset and for every wheel of the vehicle. Thus, both routines must be called for every single constraint. For the sake of brevity, the pseudocode algorithms included in this work describe only the most important details of the routine's implementation.

Both subroutines share the same interface: the function receives as input arguments the values, in the current configuration, of the degrees of freedom of all nodes involved, their initial positions, and some constant values; as output, the routine expects the newly computed values for the degrees of freedom of the slave node and the derivatives of the constraint with respect to all the degrees of freedom involved (corresponding to both the slave and master nodes). More precisely, the input variables involved in both subroutines are:
\begin{description}
  \item[\texttt{time}:] simulation time of the current iteration.
  \item[\texttt{velocity}:] the velocity $v$ of the vehicle (assumed constant for the full simulation).
  \item[\texttt{s0}:] initial longitudinal coordinates ($s$) for all the nodes involved.
    It is a one-dimensional array with length \texttt{n}\footnote{The Abaqus interface for the MPC subroutine
      does not include \texttt{velocity} and \texttt{s0} as variables: they have been explicitly included here as input variables for the sake of simplicity.
      However, they have been passed to the interface, packed into the \texttt{field} array provided by Abaqus (see the Abaqus documentation \cite{Abaqus610} for further details)}.
  \item[\texttt{X}:] initial positions and orientations for all the nodes.
    It is a two-dimensional array with dimensions $6\times\texttt{n}$ (the first three rows correspond to displacements, the last three to rotations/orientations).
  \item[\texttt{u}:] degrees of freedom values for all the nodes. It is a two-dimensional array and has the same dimensions and array distribution as \texttt{X}.
  \item[\texttt{n}:] total number of nodes involved in the constraint.
\end{description}
The output variables expected by the routine are:
\begin{description}
  \item[\texttt{ue}:] the newly computed degrees of freedom values for the slave node. It is a one-dimensional array with length 6 (three displacements and three rotations).
  \item[\texttt{A}:] the derivatives of the slave degrees of freedom with respect to the degrees of freedom involved. It is three-dimensional array with dimensions $6\times6\times\texttt{n}$, such as \texttt{A(i,j,k)} is the derivative of the \texttt{i}-th degree of freedom of the slave node respect to the \texttt{j}-th degree of freedom of the \texttt{k}-th node (it could be slave or master).
\end{description}

\algrenewcommand\algorithmicprocedure{\textbf{subroutine}}
\begin{algorithm}
\caption{Abaqus MPC for the restriction defined in Section \ref{sec:constraint1}} \label{alg:constraint_1}
\begin{algorithmic}[1]
  \For{every vehicle wheelset}
\Procedure {Mpc}{$\texttt{time}$, $\texttt{velocity}$, $\texttt{s0}$, $\texttt{X}$, $\texttt{u}$, $\texttt{n}$, $\texttt{ue}$, $\texttt{A}$}

  \State $\texttt{sm}=\texttt{s0(1)}+\texttt{velocity}\ast\texttt{time}$
  \State Find \texttt{n1} and \texttt{n2} such as $\texttt{s1}=\texttt{s0(n1)} \leq \texttt{sm} < \texttt{s2}=\texttt{s0(n2)}$
  \State Find \texttt{m1} and \texttt{m2} such as $\texttt{s0(m1)} \leq \texttt{s0(1)} < \texttt{s0(m2)}$
  \State $\bm{X}_m=\texttt{X(1:3,1)}$, $\bm{X}_1=\texttt{X(1:3,n1)}$, $\bm{X}_2=\texttt{X(1:3,n2)}$
  \State $\bm{\Phi}_m=\texttt{X(4:6,1)}$, $\bm{\Phi}^0_1=\texttt{X(4:6,m1)}$, $\bm{\Phi}^0_2=\texttt{X(4:6,m2)}$
  \State $\bm{u}_1=\texttt{u(1:3,n1)}$, $\bm{u}_2=\texttt{u(1:3,n2)}$
  \State $\bm{\phi}_1=\texttt{u(4:6,n1)}$, $\bm{\phi}_2=\texttt{u(4:6,n2)}$

  \State $\zeta = (2\ast\texttt{sm} - \texttt{s1} - \texttt{s2})/(\texttt{s1} - \texttt{s2})$, $N_1 = N_1(\zeta)$ and $N_2 = N_2(\zeta)$

  \State Computation of $\bm{E}^b_i$ ($i=\lbrace x,y,z\rbrace$) by means of $\bm{\Phi}^0_1$ and $\bm{\Phi}^0_2$
  \State Computation of $c^m_y$, $c^m_z$ and $\bm{\phi}_m$ according to \eqref{ch3:eq:excen} and \eqref{eq:interp_rotation}
  \State Computation of $\bm{\rho}_m$ according to \eqref{ch3:eq:Rhom} and \eqref{ch3:eq:rest_disp1_bis}
  \State \texttt{ue(1:3)}$=N_1\,\left(\bm{X}_1 + \bm{u}_1\right)+ N_2\,\left(\bm{X}_2 + \bm{u}_2\right) + \bm{\rho}_m - \bm{X}_m$
  \State \texttt{ue(4:6)}$=\bm{\phi}_m=N_1\,\bm{\phi}_1 + N_2\,\bm{\phi}_2$
  \State $\texttt{A(1:3,1:3,1)}=\bm{1},\quad \texttt{A(1:3,4:6,1)}=\bm{0}$ 
  \State $\texttt{A(4:6,1:3,1)}=\bm{0},\quad \texttt{A(4:6,4:6,1)}=\bm{1}$

  \State $\texttt{A(1:3,1:3,n1)}=-N_1\bm{1},\quad\texttt{A(1:3,4:6,n1)}=N_1\widehat{\bm{\rho}_m} \bm{H}(\bm{\phi}_m) \bm{T}(\bm{\phi}_1)$
  \State $\texttt{A(4:6,1:3,n1)}=\bm{0},\quad\texttt{A(4:6,4:6,n1)}=-N_1 \bm{H}(\bm{\phi}_m) \bm{T}(\bm{\phi}_1)$

  \State $\texttt{A(1:3,1:3,n2)}=-N_2\bm{1},\quad\texttt{A(1:3,4:6,n2)}=N_2\widehat{\bm{\rho}_m} \bm{H}(\bm{\phi}_m) \bm{T}(\bm{\phi}_2)$
  \State $\texttt{A(4:6,1:3,n2)}=\bm{0},\quad\texttt{A(4:6,4:6,n2)}=-N_2 \bm{H}(\bm{\phi}_m) \bm{T}(\bm{\phi}_2)$

  \For{\texttt{i} = \texttt{2 : n}}
    \If{\texttt{i} $\neq$ \texttt{n1} and \texttt{i} $\neq$ \texttt{n2}}
      \State $\texttt{A(1:6,1:6,i)}=\bm{0}$
    \EndIf
  \EndFor

\EndProcedure
  \EndFor
\end{algorithmic}
\end{algorithm}

\begin{algorithm}
\caption{Abaqus MPC for the restriction defined in Section \ref{sec:constraint2}} \label{alg:constraint_2}
\begin{algorithmic}[1]
  \For{every vehicle wheelset and for $r_1$ and $r_2$}
\Procedure {MPC}{$\texttt{time}$, $\texttt{velocity}$, $\texttt{s0}$, $\texttt{X}$, $\texttt{u}$, $\texttt{n}$, $\texttt{ue}$, $\texttt{A}$}

  \State $\texttt{sm}=\texttt{s0(1)}+\texttt{velocity}\ast\texttt{time}$
  \State Compute irregularities $I_y(s)$, $I_z(s)$, $I_\theta(s)$ from a pre-loaded table.
  
  \State $\bm{X}_r=\texttt{X(1:3,1)}$, $\bm{X}_m=\texttt{X(1:3,2)}$, $\bm{\Phi}_r=\texttt{X(4:6,2)}$, $\bm{\Phi}_m=\texttt{X(4:6,2)}$
  \State $\bm{u}_m=\texttt{u(1:3,2)}$, $\bm{\phi}_m=\texttt{u(4:6,2)}$
  \State Computation of $c^r_y$, $c^r_z$ and $\bm{\rho}_r$ according to \eqref{eq:offsets_rails} and \eqref{eq:rest_disp5bis}
  \State \texttt{ue(1:3)}$=\bm{X}_m + \bm{u}_m + \bm{\rho}_r -\bm{X}_r$
  \State \texttt{ue(4:6)}$=\texttt{spurrier}(\bm{\Lambda}(\bm{\phi}_m)\,\bm{\Lambda}(\Gamma_\theta\bm{e}^r_x))$

  \State $\texttt{A(1:3,1:3,1)}=\bm{1},\quad \texttt{A(1:3,4:6,1)}=\bm{0}$ 
  \State $\texttt{A(4:6,1:3,1)}=\bm{0},\quad \texttt{A(4:6,4:6,1)}=\bm{1}$

  \State $\texttt{A(1:3,1:3,2)}=-\bm{1},\quad \texttt{A(1:3,4:6,2)}=\widehat{\bm{\rho}_r}$ 
  \State $\texttt{A(4:6,1:3,2)}=\bm{0},\quad \texttt{A(4:6,4:6,2)}=-\bm{1}$

\EndProcedure
  \EndFor
\end{algorithmic}
\end{algorithm}

Finally, it is important to remark that, as can be appreciated in Algorithm \ref{alg:constraint_2}, the
track irregularities at a certain section $s$ are computed by interpolating in the pre-loaded lookup table.
It contains the values of $I_y$, $I_z$, and $I_\theta$ for both rails referred to a set of longitudinal coordinates. This table is stored in a comma-separated value format in an external file and pre-loaded
at the beginning of the simulation by properly implementing the \texttt{Uexternadb} subroutine provided by Abaqus
(see \cite{Abaqus610} for further details).

\section*{Acknowledgements}
Pablo Antolin acknowledges the financial support of the Swiss National Science Foundation through the project FLAS$_h$ with no. 200021\_214987.

Khanh Nguyen and Jose M.\ Goicolea are grateful for the funding from the Europe's Rail Joint Undertaking under Horizon Europe research and innovation programme under grant agreement No.\ 101121765 (HORIZON-ER-JU-2022-ExplR-02). Views and opinions expressed are, however, those of the author(s) only and do not necessarily reflect those of the European Union or Europe’s Rail Joint Undertaking. Neither the European Union nor the granting authority can be held responsible for them.

\bibliography{mybibfile}

\end{document}